\documentclass{amsart}
\usepackage[utf8]{inputenc}
\usepackage[margin=.5in,left=.5in]{geometry}
\usepackage{setspace}
\usepackage{graphicx}
\usepackage{amssymb, amsmath, amsthm, graphics}
\usepackage{mathabx,epsfig}
\usepackage[mathscr]{euscript}
\usepackage{diagbox}
\usepackage{tikz}
\usetikzlibrary{calc}
\usetikzlibrary{matrix}
\usepackage[T1]{fontenc}
\usepackage[sc,osf]{mathpazo}
\usepackage{times}
\usepackage{a4wide}  
\usepackage{slashed}
\usepackage{latexsym,amsthm,amsfonts,amsmath,mathrsfs,amssymb}
\usepackage{booktabs} 
\usepackage[unicode,implicit]{hyperref}

\usepackage{latexsym}
\usepackage[all]{xy}
\usepackage{color}

\usepackage{bbold}

\usepackage{tikz}
\xyoption{all}
\usepackage{euscript}
\usepackage{multirow}
\usepackage{etex, pictexwd,dcpic}
\usetikzlibrary{positioning}
\usetikzlibrary{shapes.geometric}
\usetikzlibrary{shapes.misc}
\usetikzlibrary{calc}
\usetikzlibrary{positioning}

 \usepackage{enumerate}

\usepackage{pgflibraryarrows}
\usepackage{pgflibrarysnakes}

\usetikzlibrary{trees} 
\usetikzlibrary[trees] 

\usepackage{stmaryrd}    

\def\acts{\mathrel{\reflectbox{$\righttoleftarrow$}}}

\newtheorem{theorem}{Theorem}
\newtheorem{definition}[theorem]{Definition}
\newtheorem{lemma}[theorem]{Lemma}
\newtheorem{corollary}[theorem]{Corollary}
\newtheorem{proposition}[theorem]{Proposition}

\newtheorem{remark}{Remark}

\newtheorem{example}{Example}

\numberwithin{equation}{section}

\usepackage{chngcntr}
\counterwithout{table}{section}

\renewcommand{\(}{\begin{equation*}}
\renewcommand{\)}{\end{equation*}}
\newcommand{\bea}{\begin{eqnarray*}}
\newcommand{\eea}{\end{eqnarray*}}

\def\endofproof {\hfill{$\Box$}\\}

\newcommand{\cL}{\ensuremath{\mathcal L}}

\def\S{\ensuremath{\ES{S}}}

\def\S{\ensuremath{\ES{S}}}

\newcommand{\beq}{\begin{equation}}
\newcommand{\eeq}{\end{equation}}

\newcommand{\onto}{\twoheadrightarrow}

\newcommand{\into}{\hookrightarrow}

\newcommand{\op}[1]{\ensuremath{\operatorname{#1}}}

\newcommand{\ES}[1]{\ensuremath{\EuScript{#1}}}

\newcommand{\theproof}{\noindent {\bf Proof.\ }}

\numberwithin{equation}{section}

\usepackage{cleveref}

\crefformat{section}{\S#2#1#3} 
\crefformat{subsection}{\S#2#1#3}
\crefformat{subsubsection}{\S#2#1#3}

\renewcommand{\(}{\begin{equation}}
\renewcommand{\)}{\end{equation}}

\def\1{{\bf 1}}

\def\<{\langle}
\def\>{\rangle}

\numberwithin{equation}{section}


\newcommand{\R}{\ensuremath{\mathbb R}}
\newcommand{\RR}{\ensuremath{\mathbb R}}

\newcommand{\ZZ}{\ensuremath{\mathbb Z}}
\newcommand{\Z}{\ensuremath{\mathbb Z}} 
\newcommand{\QQ}{\ensuremath{\mathbb Q}}

\newcommand{\BB}{\ensuremath{\mathbf B}}

\newcommand{\sh}{\ensuremath{\mathcal{S}\mathrm{h}}}

\newcommand{\cartsp}{\mathcal{C}\mathrm{art}\mathcal{S}\mathrm{p}}
\newcommand{\map}{\mathrm{Map}}

\begin{document}

\title{
Twisted differential KO-theory}


\author{Daniel Grady}
\address{Mathematics Program, Division of Science, New York University Abu Dhabi,
 Saadiyat Island, Abu Dhabi, UAE 
}
\email{djg9@nyu.edu}
\author{Hisham Sati}
\address{Mathematics Program, Division of Science, New York University Abu Dhabi,
 Saadiyat Island, Abu Dhabi, UAE }
\email{hsati@nyu.edu}

\maketitle

\begin{abstract} 
We provide a systematic approach to twisting differential  KO-theory leading to a construction 
of the corresponding twisted differential Atiyah-Hirzebruch spectral sequence (AHSS). We relate 
and contrast the degree two and the degree one twists, whose description involves appropriate 
local systems.  Along the way,  we provide a complete and explicit  identification of the differentials  
at the $E_2$ and $E_3$ pages in the topological case, which  has been missing in the literature and 
which is  needed for the general case. 
 The corresponding differentials in the refined theory  reveal   an intricate interplay between topological 
 and geometric data, the former involving the flat part and the latter requiring  the construction of the 
 twisted differential Pontrjagin character. 
 We   illustrate with examples and applications from geometry, topology and physics. 
  For instance, quantization conditions show  how to lift  differential $4k$-forms to twisted differential KO-theory 
  leading to integrality results, while  considerations of anomalies in type I string theory allow for 
 characterization of twisted differential Spin structures. 
  \end{abstract}

\tableofcontents


\section{Introduction}

\medskip
In a previous work \cite{GS-KO} we provided a comprehensive account of differential 
KO-theory $\widehat{\rm KO}$, the motivation for which was, to a large extent, to lay the 
groundwork and foundation for the current paper. Indeed, here we build on \cite{GS-KO} in order to 
explicitly construct twists for $\widehat{\rm KO}$ and provide computational techniques, which 
extend the literature even in the purely topological case.

\medskip
Twisted KO-theory has been originally considered by Donovan and Karoubi in \cite{DK}, 
with  twists  being in bijective correspondence with the set of elements in 
\(
\label{3twists} 
\ZZ/8 \times H^1(X; \Z/2) \times H^2(X; \ZZ/2)\;.
\) 
 In \cite{Ros}, Rosenberg considers twistings of KO-theory from the point of  view of operator 
 algebras, with the twist involving a real analogue of the Dixmier-Douady class or the 
 Steifel-Whitney class, in degree two.  This gives the twisted KO-groups $\op{KO}^j _{w_2}(X)$ 
 via real continuous-trace algebras. As pointed out by Rosenberg,  these twisted groups have 
 applications in the Poincar\'e duality theorem for KO of non-Spin manifolds. 
It is  then studied by Mathai, Murray, and Stevenson in \cite{MMS} from the point of view of bundle 
gerbe modules  for the purpose of applications in string theory. A treatment from the homotopy 
 point of view is given by Antieau, Gepner, and Gomez in \cite{AGG}, where uniqueness of 
 twisted KO-theory is established. The equivalence of the definitions of twisted KO-theory 
 from an operator algebra and homotopy theory points of view is proved by Hebestreit and Sagave
 in \cite{HSa}, relating to the spectrum model for KO given by Joachim \cite{Jo}\cite{Jo2} and 
 extending  the equivalence in \cite{AGG} to the level of $E_\infty$-maps. 
 A model for the twisted KO-spectrum as well as of the twisted Spin cobordism 
spectrum  in relation to metrics of positive scalar curvature are given by Hebestreit and Joachim
in \cite{HJ}. 

\medskip
From a geometric  (and to some extent, analytic) point of view, there are several important occurrences 
of twisted  KO-theory.  First,  any closed $n$-dimensional manifold $M$ admits a fundamental class 
$[M] \in KO_n(M, {or}(M))$ in its connective KO-homology twisted by ${or}(M) \to B{\rm O}$ 
given by the first two Stiefel-Whitney classes $(w_1, w_2)$ of $M$. Twisted KO-theory also 
appears in the study of positive scalar curvature. The Gromov-Lawson-Rosenberg conjecture 
(see \cite{Ros3}\cite{RS}\cite{Sto}) asserts that the vanishing of a certain index class in the twisted 
KO-homology of the $C^*$-algebra of the fundamental group, 
${\rm ind}(M) \in KO_d(C^*(\pi_1(M)), {or}(M))$, is equivalent to the 
manifold $M$ having a metric of positive scalar curvature. Connective versions also appear. 
If the universal cover $\widetilde{M}$ of $M$ admits a Spin structure then the twist ${or}(M)$ 
factors through the classifying map $c: M \to B\pi_1(M)$. 
Stolz conjectured in \cite{Sto} that the vanishing of the class $c_*[M]\in ko_d(B\pi_1(M), {or}(M))$
 is sufficient for the existence of a metric of positive scalar curvature on $M$.
It was pointed out in \cite{HSa} that it would be desirable to calculate such twisted groups. 
 In a slightly different formulation, twisted KO-groups of classifying spaces of fundamental groups
 have been considered in this context, e.g.,  by Botvinnik, Gilkey, and Stolz in \cite{BGS}, where
  a homological version of the Atiyah-Hirzebruch spectral sequence (AHSS) for 
   twisted KO-theory is used to check the conjecture for finite groups with periodic cohomology. 
Currents twisted by orientation bundles using Stiefel-Whitney classes have been developed by Harvey and Zweck
\cite{HZ}.

\medskip
For KO-theory one twists by elements of the set \eqref{3twists}, with $H^m(X; \ZZ/2)$, $m=1,2$,
 associated to a map $f: X \to K(\Z/2, m)$.  
Composing with the inclusion map $i: K(\ZZ/2, m) \to B{\rm GL}_1{\rm KO}$
(see below) gives $f$-twisted KO-theory ${\rm KO}(X)_{i \circ f}$.
There is a decomposition of infinite loop spaces \cite{MST}\cite{AGG}
$$
B{\rm GL}_1{\rm KO} \simeq K(\Z/2, 1) \times K(\Z/2, 2) \times BB{\rm SO}_\otimes\;.
$$
 In \cite{AGG} it was shown that there is a natural isomorphism
of groups 
$$
[K(\ZZ/2, 2), B{\rm GL}_1{\rm KO}] \cong [K(\ZZ/2, 2), K(\ZZ/2, 2)] \cong \Z/2\;.
$$
This means that any two maps from $K(\Z/2, 2)$ to $B{\rm GL}_1{\rm KO}$ differ by an endomorphism 
of $K(\Z/2, 2)$, up to homotopy. Furthermore, it was shown in \cite{AGG} that 
any map $j: K(\ZZ/2, 2) \to B{\rm GL}_1{\rm KO}$ is either homotopically 
trivial or is equivalent to the canonical inclusion. This then establishes
 unique and nontrivial definition of twisted KO-theory.

 \medskip
 One could  start with twisted KO-theory and then refine it differentially or,  alternatively, start with 
 differential KO-theory $\widehat{\rm KO}$ from \cite{GS-KO} and then twist it. 
 Note that whether we call a theory ``differential twisted" or ``twisted differential" 
does not matter, as twistings of differential theories are differential twists of the
underlying topological theory (see also \cite{GS-KO}).
In the approach of \cite{GS5}, for twisted KO-theory ${\rm KO}_{{\rm tw}}$, 
differential KO-theory $\widehat{\rm KO}$ and twisted differential KO-theory 
$\widehat{\rm KO}_{\widehat{\rm tw}}$, we have the following schematic picture   
$$
\xymatrix{
\widehat{\rm KO}_{\widehat{\rm tw}} \ar@{~>}[d]_{\widehat {\rm tw}=0} \ar[rrr]^{I_I} 
 &&& {\rm KO}_{\rm tw} \ar@{~>}[d]^{{\rm tw}=0} \\
 \widehat{\rm KO} \ar[rrr]^I &&&
 {\rm KO}
}
$$
One of our goals is to provide computational tools for twisted differential KO-theory,
which allow for explicit computations  for various smooth manifolds. We have established an 
Atiyah-Hirzebruch spectral sequence (AHSS) for differential KO-theory in \cite{GS-KO}. 
As in the complex case \cite{GS5}, we will have a correspondence diagram of transformations
 of the corresponding spectral sequences 
$$
\xymatrix{
 \widehat{\rm AHSS}_{\widehat {\rm tw}} \ar@{~>}[d]_{\widehat {\rm tw}=0} \ar[rrr]^{I_I} 
 &&& {\rm AHSS}_{\rm tw} \ar@{~>}[d]^{{\rm tw}=0} \\
 \widehat{\rm AHSS} \ar[rrr]^I &&&
 {\rm AHSS}
}
$$
The AHSS for twisted KO-theory ${\rm AHSS}_{\rm tw}$ has been set up in \cite{Ros}  
and used briefly  in \cite{MMS} and by Doran, Mendez-Diez, and Rosenberg in\cite{DMDR} to 
calculate 2-dimensional examples. When the twist is turned off we reproduce the results in 
\cite{GS-KO} for $\widehat{\rm AHSS}$, when the geometry is turned off we reproduce the 
results in \cite{Ros}, while when both the twist and the differential refinement are turned off, 
we reproduce the classical construction of Atiyah and Hirzebruch \cite{AH}.

   \medskip
 In our quest to explicitly work out the differentials in $\widehat{\rm AHSS}_{\widehat {\rm tw}}$,
 we realized that the corresponding differentials in the twisted topological case ${\rm AHSS}_{{\rm tw}}$    
 had not been fully  worked out, even in the lowest pages, namely the $E_2$ and  $E_3$ pages. 
Hence we explicitly work out these differentials, which turns out to be quite delicate. 
In the absence of a K\"unneth spectral sequence  as, for instance, the calculations of
 Dobson in \cite{Dob} highlight,  
we had to do many of the  calculations from scratch using real projective spaces and Klein bottles. 
An excellent exposition of KO-theory of (stunted) real projective spaces is given by Dugger in \cite[Sec. 32]{Dug2}
and we have found the twisted Thom isomorphism in KO-theory established by Hebestreit and Joachim in
\cite[Prop. 4.1.7]{HJ} to be useful in  computing the twisted $\op{KO}$-theory of $\RR P^n$.
The cohomology groups of the Klein bottle over $\Z$ and $\Z/2$ are calculated by Palermo in \cite[p. 185]{Pa}. 
For higher-dimensional Klein bottles $K_n$, not much seems to have been known previously in the literature. 
For instance,  for $K_4$, the $\Z/2$-cohomology  and KO-theory are briefly indicated via the AHSS by 
Davis in \cite{Dav}, where the author  expresses doubts on the utility of such results. We hope that our perspective 
a posteriori provides  such a  motivation. To highlight the utility, we provide an inductive argument on the dimension 
$n$ of the Klein bottle $K_n$. This generalizes the computations 
in \cite{DMDR} for the circle, torus, and M\"obius strip. In contrast, manifolds  satisfying $w_1^2 \neq 0$ 
are described by Lusztig, Milnor, and Peterson in \cite{LMP} in their work on semi-characteristics.

 \medskip
The group of twists is a non-split abelian extension of $H^2(X; \ZZ/2)$ by $H^1(X; \ZZ/2)$ \cite{DK}
$$
\xymatrix{
0 \longrightarrow H^2(X; \ZZ/2) \longrightarrow {\rm GBR}(X) 
\longrightarrow \ZZ/8 \oplus H^1(X; \Z/2) \longrightarrow 0 
}.
$$ 
The graded real Brauer group ${\rm GBR}(X)$ has been determined in many cases, e.g.,
for the real projective plane ${\rm GBR}(\R P^2) \cong \Z/8 \oplus \Z/4$  by Karoubi, Schlichting, and Weibel 
\cite[Ex. 2.5]{KSW}\cite[Ex. 2.8]{KW2}  using techniques developed therein.
Distler, Freed and Moore \cite{DFM2} \cite{DFM1} consider both periodic KO-theory and its 
connective cover $\op{ko}$, as well as a variant multiplicative cohomology theory $R^{-1}$, which 
is the Postnikov section $\op{ko}\langle 0, \cdots, 4\rangle$, describing the set of twists for complex
K-theory.  It was indicated that twistings and differential refinements (separately) of the latter  
are inherited from those of $\op{ko}$ or KO. Indeed, we will show how to combine the two aspects within
the twisted differential theory. On the other hand, it was asserted (with proof to appear) in \cite{DFM2} that 
$R^{w_0 -1}(B\Z/2; \Z)$ can be interpreted as the group of universal twistings of KO-theory modulo Bott
periodicity, the latter being the super Brauer group.

 \medskip
 Twists of KO-theory in relation to applications to orbifolds and orientifolds 
from string theory is  treated by Freed in \cite{Fr-lec}.  In this context, Atiyah's KR-theory \cite{Real}
which is, in some sense, an intermediate ground between KO-theory and K-theory, is also 
important \cite{DFM2}\cite{DFM1}\cite{DMDR}\cite{DMDR2}. Twists for KR-theory,  have been 
developed generally via groupoids by Moutuou in \cite{Mou1}\cite{Mou2}\cite{Mou3}.
Twisted KR-homology aimed towards geometry and physics applications, is described by
Hekmati,  Murray,  Szabo, and Vozzo in \cite{HMSV}. 
However, as in \cite{GS-KO}, our focus in this paper will be on KO-theory itself as a cohomology theory. 

\medskip
We also highlight that having an axiomatic 
approach to twisted differential cohomology -- as developed by Bunke and Nikolaus \cite{BN} -- does 
not guarantee ability to  explicitly construct such a theory, as we highlighted and demonstrated 
in \cite{GS4}\cite{GS5}\cite{GS6}. This is the case even in the purely topological case where 
the constructions can be highly nontrivial; see \cite{ABG}\cite{SW}\cite{LSW} for illustrations. 
We use moduli stacks of bundles with connection representing differential cohomology
\cite{Cech}\cite{Urs}\cite{SSS3} \cite{FSS1}\cite{FSS2} adapted and extended to $\Z/2$ coefficients
\cite{GS4}\cite{GS6}. 

\medskip
The paper is organized as follows. 
In Sec. \ref{Sec-tw}, we provide detailed construction of twisting differential KO-theory. In more detail, 
in Sec. \ref{Sec-twdiff} we describe the the twists of differential KO-theory, then in Sec. \ref{Sec-tPh} 
we construct the twisted KO-theory Pontrjagin character via descent. Having constructed twisted differential 
KO-theory, we provide its properties in Sec. \ref{Sec-prop}. This then allows us to construct the twisted 
differential Pontrjagin character in Sec. \ref{Sec-tdPh}. Then we move to the computational aspects in the 
form of the Atiyah-Hirzebruch spectral sequence (AHSS) for twisted differential KO-theory in Sec. \ref{Sec-AHSS}. 
We start with the AHSS for twisted topological KO-theory in Sec. \ref{Sec-tAHSS}, where we establish the explicit form of 
the differentials, which allows us to describe the ones corresponding to the AHSS for twisted differential KO-theory 
in Sec. \ref{Sec-AHSS-hatKO}. We end in Sec. \ref{Sec-apps} with applications to geometry, topology, and 
physics. It is interesting that considering field  quantization in type I string theory via the AHSS leads to interesting 
topological and geometric consequences, including a splitting for low-dimensional manifolds, lifting forms to 
$\widehat{\rm KO}_{\widehat{\rm tw}}$, integrality of such forms leading to 
directly obtaining Rokhlin's theorem on the divisibility of the 
signature of 4-manifolds \cite{Rok}, and characterization of anomalies in string theory 
in relation to twisted differential Spin structures.

\section{Twisting differential KO-theory} 
\label{Sec-tw}

\subsection{The twisted KO-theory Pontrjagin character via descent} 
\label{Sec-tPh}

In this section we define the twisted Pontrjagin character for twisted KO-theory using descent, motivated 
by the model for twisted cohomology developed  by \cite{ABGHR}. The reason for this perspective is twofold. 
First, it prepares us for the differential refinement. Second, it circumvents much of the difficulty one encounters 
in choosing models for the theory (see, e.g., \cite{CMW}).

\medskip
Instead of considering the twisted $\op{KO}$-theory spectrum as a global object, we consider the spectrum 
locally, over its space of twists, as we did in \cite{GS5}. This allows for bringing in ideas familiar from 
geometry of usual bundles.  Abstractly, twisted spectra are elements in the tangent $\infty$-topos 
$T(\infty\text{-}\mathcal{G}{\rm pd})$, which is itself an $\infty$-topos (see \cite{Joy}\cite{Urs}\cite{BNV}\cite{HNP}).    
From the Giraud-Lurie-Rezk axioms \cite{Lur}\cite{Rez} for $\infty$-toposes, one can use the principle of descent to glue 
bundles of spectra from local data. The most important incarnation of this principle for us will be the following.

\begin{proposition}[Bundles of spectra from local data]
\label{prop-00top}
Let $M$ be a topological manifold and let $\{U_{\alpha}\}$ be a good open cover of $M$. Let $\mathcal{R}$ be 
a ring spectrum with space of units $B{\rm GL}_1(\mathcal{R})$. Since each $U_{\alpha}$ is contractible, 
we have iterated pullback squares in the  tangent $\infty$-topos
\(\label{tanginftpdg}
\xymatrix@R=1.5em{
\hdots \ar@<-.1cm>[r]\ar[r]\ar@<.1cm>[r] &
\coprod_{\alpha\beta}U_{\alpha\beta}\times \mathcal{R}\ar@<.05cm>[r]\ar@<-.05cm>[r]\ar[d] &
\coprod_{\alpha}U_{\alpha}\times \mathcal{R} \ar[d] \ar[r] &  \mathcal{R}_{\sigma}\ar[r]\ar[d] & 
\mathcal{R}\sslash {\rm GL}_1(\mathcal{R})\ar[d]
\\
\hdots \ar@<-.1cm>[r]\ar[r]\ar@<.1cm>[r] & 
\coprod_{\alpha\beta}U_{\alpha\beta}\ar@<.05cm>[r]\ar@<-.05cm>[r] &
 \coprod_{\alpha}U_{\alpha}\ar[r] & M \ar[r]^-{\sigma} &B{\rm GL}_1(\mathcal{R})
}
\)
In this case, the diagram 
$$
\left\{\vcenter{\xymatrix{
\hdots \ar@<-.1cm>[r]\ar[r]\ar@<.1cm>[r] &
\coprod_{\alpha\beta}U_{\alpha\beta}\times \mathcal{R}\ar@<.05cm>[r]\ar@<-.05cm>[r]&
\coprod_{\alpha}U_{\alpha}\times \mathcal{R}  \ar[r] &  \mathcal{R}_{\sigma}
}}\right\}
$$
is also colimiting.
\end{proposition}
\theproof
From the Nerve Theorem (see \cite[Cor. 4G.3]{Ha}), the diagram 
$$
\big\{\vcenter{\xymatrix{
\hdots \ar@<-.1cm>[r]\ar[r]\ar@<.1cm>[r] &
\coprod_{\alpha\beta}U_{\alpha\beta}\ar@<.05cm>[r]\ar@<-.05cm>[r]&
\coprod_{\alpha}U_{\alpha} \ar[r] &  M
}}\big\}
$$
is colimiting in $\mathcal{T}{\rm op}$. The ``zero section" functor 
$s:\infty\text{-}\mathcal{G}{\rm pd}\longrightarrow T(\infty\text{-}\mathcal{G}{\rm pd})$ is both 
a left and right adjoint and, in particular, preserves colimits. Hence the bottom diagram in \eqref{tanginftpdg} 
is colimiting in the tangent topos. The claim then is an immediate corollary of the axiom of descent 
(more precisely, that all colimits are Van-Kampen)  for the tangent topos. 
\endofproof

Proposition \ref{prop-00top} is a higher homotopical manifestation of a classical fact from geometry, namely, 
that locally trivial  bundles are obtained by gluing local trivializing patches. See \cite{GS4}\cite{GS5}\cite{GS6} for 
explicit illustrations for other spectra. 
We now utilize this machinery to define the twisted Pontrjagin character for twisted KO-theory. 
We will use the Eilenberg-MacLane functor $H$ to 
describe the rationalization. 

\begin{proposition}[Twisted Pontrjagin character]
Let $M$ be a topological manifold and let $\sigma:M\to B{\rm GL}_1(\op{KO})$ be a twist. There exists a bundle map
$$
{\rm Ph}_{\sigma}:\op{KO}_{\sigma}\longrightarrow H(\pi_*(\op{KO})\otimes \RR)_{\sigma^*\widetilde{{\rm Ph}}}\;,
$$
called the \emph{twisted Pontrjagin character}, 
which locally restricts to the untwisted Pontrjagin character.
\end{proposition}

\theproof
Fix a good open cover $\{U_{\alpha}\}$ of $M$. Pulling back the twisted spectrum 
$\op{KO}_{\sigma}\to M$ by the inclusion $i_{\alpha}:U_{\alpha}\into M$, gives a twisted spectrum 
$\op{KO}_{i^*_{\alpha}\sigma}$ over $U_{\alpha}$. Since the set $U_{\alpha}$ is contractible, we have an 
equivalence $f_{\alpha}: \op{KO}_{i^*_{\alpha}\sigma}\simeq \op{KO}$ over $U_{\alpha}$. Similarly, 
we have local trivializations of Eilenberg-MacLane spectra 
$g_{\alpha}:H(\pi_*({\rm KO})\otimes \RR)_{i^*_{\alpha}\sigma^*\widetilde{{\rm Ph}}}
\simeq H(\pi_*({\rm KO})\otimes \RR)$ 
over each $U_{\alpha}$. We can locally define the twisted Pontrjagin character via the local trivializations as
$$
\xymatrix@R=2em{
\op{KO}_{i^*_{\alpha}\sigma}\ar[d]_{{\rm Ph}_{i^*_{\alpha}\sigma}}\ar[rr]_-{\simeq}^{f_{\alpha}} && 
\op{KO}\ar[d]^-{\rm Ph}
\\
H(\pi_*({\rm KO})\otimes \RR)_{i^*_{\alpha}\sigma^*{\rm Ph}}&& 
\ar[ll]^-{\simeq}_-{g^{-1}_{\alpha}} H(\pi_*({\rm KO})\otimes \RR)\;.
}
$$
Moreover, since ${\rm Ph}$ is a morphism of ring spectra, it follows that the local twisted Pontrjagin character maps 
${\rm Ph}_{i^*_{\alpha}\sigma}$ are compatible   with the transition data for the bundle of spectra. By descent,
 it follows  that the local morphisms 
${\rm Ph}_{i_{\alpha}^*\sigma}:\op{KO}_{i^*_{\alpha}\sigma}\to 
H(\pi_*({\rm KO})\otimes \RR)_{i^*_{\alpha}\sigma^*{\rm Ph}}$
 induce a global morphism which is unique up to a contractible choice. 
 \endofproof
 
 We will consider the  differential Pontrjagin character in Prop. \ref{prop-diffpontch}
 and the twisted differential Pontrjagin character in Theorem \ref{thm-twdifPo}.

 \medskip
The following shows that for twists of the form $\sigma:M\to K(\ZZ/2,1)\times K(\ZZ/2,2)\into B{\rm GL}_1(\op{KO})$, 
only the degree 1 twist affects the image of the Pontrjagin character. We will highlight further effects of this 
in the following sections. Note that we deal with the reals $\RR$ rather than the rationals $\QQ$
 as we are mainly interested in bringing in geometry later. 

\begin{proposition}[Rational twists of KO-theory of spaces] \label{prop-rattw}
We have an identification 
$$
B{\rm GL}_1(\pi_*(\op{KO})\otimes \RR)\simeq K(\RR^{\times},1)\times \displaystyle{\prod_{k>0}} K(\RR,4k+1)\;.
$$
\end{proposition}
\theproof
For a ring spectrum $\mathcal{R}$, recall the model for the units \cite{MQRT} given by the strict 
pullback of $\infty$-groupoids (which models the homotopy pullback) 
$$
\xymatrix@R=1.3em{
B{\rm GL}_1(\mathcal{R})\ar[r]\ar[d] & \Omega^{\infty}\mathcal{R}\ar[d]
\\
\pi_0(\mathcal{R})^\times\ar[r] & \pi_0(\mathcal{R})\;.
}
$$
We will compute this pullback. By Bott periodicity, the graded algebra $\pi_*(\op{KO})\otimes \RR$ is 
isomorphic to  $\RR[\alpha,\alpha^{-1}]$, where $\vert \alpha \vert=4$. Since the Dold-Kan functor 
${\rm DK}:\mathcal{C}{\rm h}_+\to \infty\text{-}\mathcal{G}{\rm pd}=\mathcal{K}{\rm an}$
is (strict) right adjoint, it preserves (strict) pullbacks and we can compute
$$
{\rm GL}_1(\pi_*(\op{KO})\otimes \RR)\simeq {\rm DK}\Big(\RR^{\times}\times   \prod_{k>0}\RR[4k]\Big)\simeq
 \RR^{\times}\times \prod_{k>0}K(\RR,4k)\;.
 $$
Delooping proves the claim.
\endofproof

We have indicated in the Introduction that differential refinement commutes with twisting. 
Now we compare twisting rationalized KO-theory to rationalizing 
twisted KO-theory. 

\begin{proposition}[Rationalizing commutes with twisting] \label{Prop-commute}
Let $M$ be a topological manifold and let $\sigma=(\sigma_1,\sigma_2):M\to K(\ZZ/2,1)\times K(\ZZ/2,2)$ be a 
twist for $\op{KO}$.  
\item {\bf (i)} Then the twist $\sigma_2^*{\rm Ph}:M\to B{\rm GL}_1(\op{KO})\to 
B{\rm GL}_1(H(\pi_*(\op{KO})\otimes \RR))$ is nullhomotopic and a choice of nullhomotopy produces an equivalence
$$
H(\pi_*(\op{KO})\otimes \RR))_{\sigma^*{\rm Ph}}
\xymatrix{\ar[r]^{\simeq} &}
H(\pi_*(\op{KO})_{\sigma_1}\otimes \RR)\;,
$$
where $\pi_*(\op{KO})_{\sigma_1}$ is the local coefficient system with fiber $\pi_*(\op{KO})$, defined by $\sigma_1$. 
\item {\bf (ii)} Hence, the twisted Pontrjagin character defines a map 
$$
{\rm Ph}_{\sigma}:\op{KO}_{\sigma}\longrightarrow H(\pi_*(\op{KO})_{\sigma_1}\otimes \RR)\;,
$$
where the right hand side is
 the Eilenberg-MacLane spectrum associated to the graded local coefficient system.
\end{proposition}
\theproof
From Prop. \ref{prop-rattw}, we see that on the space of units, the map induced by ${\rm Ph}$ must send the factor 
$\sigma_2$ to a nullhomotopic map in $B{\rm GL}_1(H(\pi_*(\op{KO}\otimes \RR))$. On the other hand, the restriction 
of ${\rm Ph}$ to the factor $K(\ZZ/2,1)$ is the map $i:K(\ZZ/2,1)\to K(\RR^{\times},1)$, induced by the inclusion 
$\ZZ/2\into \RR^{\times}$. A choice of nullhomotopy $h:\sigma_2^*{\rm Ph}\simeq \ast$ gives rise to a homotopy 
$i\times h:(\sigma_1^*{\rm Ph},\sigma_2^*{\rm Ph})\simeq \sigma_1^*{\rm Ph}$ in $B{\rm GL}_1(\op{KO})$. This 
results in an equivalence of corresponding twisted spectra 
$$
H(\pi_*(\op{KO})\otimes \RR))_{\sigma^*{\rm Ph}}
\xymatrix{\ar[r]^{\simeq} &}
H(\pi_*(\op{KO})\otimes \RR)_{\sigma_1^*{\rm Ph}}\;.
$$
 Since the twists by $\sigma_1$ correspond to the action by the units $\ZZ/2=\pi_0(\op{KO})^{\times}$, the latter 
 spectrum is the Eilenberg-MacLane spectrum on the local coefficient system defined by the twist, as claimed.
\endofproof

We will use these statements when combining the twist with the differential refinement. 

\subsection{The twists of differential KO-theory} 
\label{Sec-twdiff}

We begin with a brief account of the Hopkins-Singer type \cite{HS} sheaf of spectra which leads to the differential refinement 
$\widehat{\op{KO}}$ of KO-theory (see \cite{GS-KO} for more details). Let $\op{KO}$ denote the spectrum for 
topological  $\op{KO}$-theory and let $\delta(\op{KO})$ be the associated sheaf of spectra given by the constant 
sheaf functor  
$\delta:\mathcal{S}{\rm p}\to 
{\mathcal{S}\mathrm{h}}_{\infty}({\mathcal{C}\mathrm{art}\mathcal{S}\mathrm{p}};\mathcal{S}{\rm p})$. 
Applying the functor 
$\delta$ to the  Pontrjagin character gives a morphism of sheaves of spectra 
$\delta({\rm Ph}):\delta(\op{KO})\to \delta(H(\pi_*(\op{KO})\otimes \RR))$. 
We have a corresponding pullback diagram 
$$
\xymatrix@C=4em{
{\rm diff}(\op{KO},{\rm Ph},\pi_*(\op{KO})\otimes \RR)\ar[r]\ar[d] & 
H(\tau_{\geq 0}\Omega^*(-;\pi_*(\op{KO})))\ar[d]
\\
\delta(\op{KO})\ar[r]^-{\iota \circ \delta({\rm Ph})} & H(\Omega^*(-;\pi_*(\op{KO})))
}
$$
where $\iota$ is induced by the inclusion $\delta(\pi_*(\op{KO})\otimes \RR)\into \Omega^*(-;\pi_*(\op{KO}))$. 
Then the sheaf of spectra ${\rm diff}(\op{KO},{\rm Ph},\pi_*(\op{KO})\otimes \RR)$ represents the differential 
refinement $\widehat{\op{KO}}$ of $\op{KO}$ \cite{GS-KO}. We now discuss the twists of this theory.

\medskip
The main goal of this section is to understand what it means to differentially refine the twists of the form 
$\sigma:M\to K(\ZZ/2,1)\times K(\ZZ/2,2)$, both from a homotopy theoretic and geometric points of view. 
Because the twists are torsion classes, the differential refinement needs to be 
interpreted with some care. In particular, we will see that the geometry does not explicitly see the degree 2 twists, 
since they are killed  under rationalization. This is in contrast to the case of complex $K$-theory, where the periodic 
de Rham complex  $\Omega^*[u,u^{-1}]$, $\vert u \vert=2$, becomes twisted by a degree three de Rham class $H$. 
The degree one twist here, which are seen by the differential refinement, are analogous to torsion twists of integral 
cohomology leading to twisted Deligne cohomology with a twist of degree one \cite{GS4}; we will see this explicitly below. 


\medskip
We begin with a general discussion on differential refinements of twists. In contrast to ordinary cohomology 
theories -- which have a \emph{space} of twists -- differential cohomology theories have a \emph{smooth stack} 
of twists. In \cite{BN}, Bunke and Nikolaus defined the smooth stack of twists 
 for a differential refinement of a ring spectrum ${\rm diff}(\mathcal{R},c,A)$, as the $\infty$-pullback 
 of smooth stacks 
\(\label{rdftwi}
\xymatrix@R=1.5em{
\widehat{\rm Tw}_{\widehat{\mathcal{R}}}\ar[rr]\ar[d] && {\rm Pic}_{\Omega^*(-;A)}\ar[d]
\\
{\rm Pic}_{\underline{\mathcal{R}}}\ar[rr] && {\rm Pic}_{H\Omega^*(-;A)}
}
\)
where ${\rm Pic}_{\Omega^*(-;A)}$ is the smooth stack of weakly locally constant, $K$-flat, invertible modules over 
$\Omega^*(-;A)$, ${\rm Pic}_{H(\Omega^*(-;A))}$ is the stack of locally constant, invertible sheaves of module spectra 
over $H\Omega^*(-;A)$, and ${\rm Pic}_{\delta(\mathcal{R})}$ is the stack of invertible module spectra over $\mathcal{R}$, 
regarded as a locally constant smooth stack (see \cite{BN}\cite{GS5}\cite{GS-KO} for detailed discussions). 
In particular, for $\op{KO}$,  we make the following definition.
\begin{definition}
[The stack of twists for $\widehat{{\rm KO}}$]\label{Def-stackotw}
The stack of twists for differential KO-theory is defined via diagram \eqref{rdftwi}
as 
\(
\label{Tw-KO}
\xymatrix@R=1.5em{
\widehat{\rm Tw}_{\widehat{\op{KO}}}\ar[rr]^-{\mathcal{R}}\ar[d]_-{\mathcal{I}} &&
 {\rm Pic}_{\Omega^*(-;\pi_*({\rm KO}))}\ar[d]
\\
{\rm Pic}_{\delta({\rm KO})}\ar[rr] && {\rm Pic}_{H(\pi_*({\rm KO})\otimes \RR)}
\;,
}
\)
where $\pi_*(\op{KO})$ is the graded ring of stable homotopy groups.
\end{definition}

\medskip
The smooth stack of twists $\widehat{\rm Tw}_{\widehat{{\rm KO}}}$ is analogous to the space \footnote{Actually, 
it is a direct generalization of the Picard $\infty$-space ${\rm Pic}_{\infty}({\rm KO})$, of which  
$B{\rm GL}_1({\rm KO})$ is the connected component of the identity.} 
$B{\rm GL}_1({\rm KO})$ in the topological case. As we have seen, ${\rm KO}$ admits twists by degree one and two
ordinary cohomology classes with $\ZZ/2$-coefficients. We would like to be able to twist $\widehat{{\rm KO}}$ by 
\emph{differential} cohomology classes, although a priori the torsion seems to obstruct this refinement. Nevertheless, we will 
see that the geometry is affected by both twists. The degree 1 twists will affect the forms in a fairly transparent way, 
while the degree 2 twists are much more subtle and are discussed in Prop. \ref{prop-geomdiff} in Sec. \ref{Sec-AHSS-hatKO}.

\begin{remark}[Twists of degree 1]
 These twists affect the coefficients of differential forms in the following way, as in \cite{GS4}. 
If $\sigma_1:M\to \BB\ZZ/2\simeq \delta(K(\ZZ/2,1))$ is a twist of degree 1 for topological $\op{KO}$, then the 
corresponding twist of the rationalization (see Prop. \ref{prop-rattw} and Prop. \ref{Prop-commute})
can be identified with the composite map
\(\label{rattwfromz2}
\xymatrix{
r(\sigma_1):\BB\ZZ/2\ar[r]^-{u} & \BB(\RR^{\times})^{\delta}\simeq \delta(K(\RR^{\times},1))
\; \ar@{^{(}->}[r] & {\rm Pic}_{\delta(H(\pi_*(\op{KO})\otimes \RR))}
 }, 
\) 
where $u$ is induced by the inclusion $\ZZ/2\into \RR^{\times}$. Let $\BB \RR^{\times}_{{\rm fl}\text{-}\nabla}$ be 
the moduli stack of real line bundles with flat connections. This can be presented, for example, by the stackification of the 
action groupoid $C^{\infty}(-;\RR^{\times})\acts \Omega^1_{\rm cl}(-)$ where, on every smooth manifold $M$, 
the function $g\in C^{\infty}(M;\RR^{\times})$ acts by gauge transformations 
$g:\mathcal{A}\mapsto \mathcal{A}+g^{-1}dg$. The Riemann-Hilbert 
correspondence gives an equivalence of smooth stacks
$$
\BB \RR^{\times}_{{\rm fl}\text{-}\nabla}\simeq \BB(\RR^{\times})^{\delta}\simeq \delta(K(\RR^{\times},1))\;.
$$
Thus, we can  take the map $r(\sigma_1)$ in \eqref{rattwfromz2} to have values either 
in $\BB \RR^{\times}_{{\rm fl}\text{-}\nabla}$ or in the geometrically discrete stack $\BB(\RR^{\times})^{\delta}$ 
classifying local systems.  The equivalence between the two perspective gives a twisted de Rham theorem, as in \cite{GS4}. 
\end{remark}

\medskip
In more detail, let $\mathcal{L}\to M$ be a real graded line bundle with fiber the rationalized coefficients 
$\RR[\alpha,\alpha^{-1}]\cong \pi_*(\op{KO})\otimes \RR$ where $\vert \alpha\vert=4$. Fix a flat connection $\nabla$ 
on $\mathcal{L}$ and consider the corresponding classifying map 
$(\mathcal{L},\nabla):M\to \BB \RR^{\times}_{{\rm fl}\text{-}\nabla}$, regarded as a twist for the rationalization 
of $\op{KO}$ (see Prop. \ref{prop-rattw} and Prop. \ref{Prop-commute}). This gives a 4-periodic twisted de Rham complex
\(\label{twfrmdg1}
\Omega^*(-;\mathcal{L}):=\xymatrix{
\Omega^0(-;\mathcal{L})\ar[r]^-{\nabla} & \Omega^1(-;\mathcal{L})\ar[r]^-{\nabla} 
& \Omega^2(-;\mathcal{L})\ar[r] & \hdots \;.
}
\)
On the other hand, we can fix an equivalence between $(\mathcal{L},\nabla)$ and a twist 
$\sigma_1^{\prime}:M\to \BB (\RR^{\times})^{\delta}\into \BB \RR^{\times}_{{\rm fl}\text{-}\nabla}$ 
classifying a local system with fiber $\pi_*(\op{KO})\otimes \RR$. This 
equivalence gives rise to a twisted de Rham equivalence 
\(\label{deReq}
\xymatrix{
d:H(\pi_*(\op{KO})_{\sigma_1}\otimes \RR)
\ar[rr]^-{\simeq}_-{\tiny \text{twisted de Rham}}&& H(\Omega^*(-;\mathcal{L}))
}.
\)
If the twist $\sigma^{\prime}_{1}$ further factors through the map $u$ in \eqref{rattwfromz2}, i.e.,
 $\sigma^{\prime}_1=u\sigma_1$, then in total we have a homotopy commutative diagram 
\(
\label{diagtwstko}
\hspace{1.5cm}
\xymatrix@C=3.7em{
& & \Omega^1_{\rm cl}\ar[dr]^{\bf 3}\ar[r]^-{\bf 2} & {\rm Pic}_{\Omega^*(-;\pi_*(\op{KO}))}\ar@/^2pc/[dd]^-{\bf 4}
\\
M\ar[r]^-{\sigma_1}_{\bf 6}\ar@/^1pc/[rru]_{\mathcal{A}}^-{\bf 1}\ar[dr] & \BB\ZZ/2\ar[r]^-{u}_{\bf 7} \ar[d]^-{\bf 9} & 
\BB(\RR^{\times})^{\delta}\ar[r]_{\bf 8}^-{\simeq}\ar[d] & \BB \RR^{\times}_{{\rm fl}\text{-}\nabla}\ar[d]_{\bf 5} & 
\\
& {\rm Pic}_{\delta(\op{KO})}\ar[r] &  {\rm Pic}_{\delta(H(\pi_*(\op{KO})\otimes \RR))}\ar[r]^{\simeq} & 
 {\rm Pic}_{\delta(H(\Omega^*(-;\pi_*(\op{KO}))))} 
}
\)
where the composite ${\bf 2}\circ {\bf 1}$ picks a globally defined flat connection $\nabla=d+\mathcal{A}$ and sends it to the 
corresponding twisted de Rham complex \eqref{twfrmdg1}. The composite ${\bf 3}\circ {\bf 1}$ is the classifying map for the line 
bundle with flat connection $(\mathcal{L},\nabla)$. 
The composite ${\bf 8}\circ {\bf 7}\circ {\bf 6}$ is the factorization of the twist $\sigma^{\prime}_1$ through the map 
$u:\BB\ZZ/2\to \BB(\RR^{\times})^{\delta}$, and the homotopy in ${\rm Pic}_{\delta(H(\Omega^*(-;\pi_*(\op{KO})))}$ 
that fills the outer diagram is identified with the equivalence \eqref{deReq}. By the universal property, this gives rise to a map 
to the pullback from the stack of twists of Def. \ref{Def-stackotw}
$$
\hat{\sigma}_1:M\longrightarrow \widehat{\rm Tw}_{\widehat{\op{KO}}}\;.
$$
Working universally, we have the following more formal statement.

\begin{proposition}[Differential degree 1 twists] 
\label{sttwdeg1}
Let $\BB(\ZZ/2)_{\nabla}$ denote the smooth stack given by the pullback 
\(\label{Bz2tw}
\xymatrix@R=1.3em@C=3em{
\BB(\ZZ/2)_{\nabla}\ar[r] \ar[d]_{\mathcal{I}^{\prime}} & \Omega^1_{\rm cl}\ar[d]
\\
\BB\ZZ/2\ar[r] & \BB\RR^{\times}_{{\rm fl}\text{-}\nabla}\;.
}
\)
\item {\bf (i)} There exists a map $\iota:\BB(\ZZ/2)_{\nabla}\to \widehat{{\rm Tw}}_{\widehat{\op{KO}}}$ which refines the 
inclusion $i:K(\ZZ/2,1)\into B{\rm GL}_1(\op{KO})$ in the sense that $i\mathcal{I}^{\prime}=\mathcal{I}\iota$, 
with $\mathcal{I}$ the canonical map in diagram \eqref{Tw-KO}. 
\item {\bf (ii)} Moreover, given a map 
$\hat{\sigma}_1:M\to \BB(\ZZ/2)_{\nabla}$, the induced map 
$$
\xymatrix{
\mathcal{R}\circ \iota\circ \hat{\sigma}_1:M \ar[r] &\BB(\ZZ/2)_{\nabla}
\ar[r] & \widehat{\rm Tw}_{\widehat{\op{KO}}} \ar[r] &
 {\rm Pic}_{\Omega^*(-;\pi_*(\op{KO}))}
 }
$$
picks out a twisted de Rham complex of the form \eqref{twfrmdg1} and the induced homotopy can be identified with an 
equivalence (twisted de Rham theorem)
$$
\xymatrix{
d:H(\Omega^*(-;\mathcal{L}))\ar[r]^-{\simeq} &  H(\pi_*(\op{KO})_{\sigma_1}\otimes \RR)
},
$$
with $\sigma_1:M\to K(\ZZ,2,1)$ the underlying topological twist.
\end{proposition}
\theproof
The proof amounts to a more careful consideration of diagram \eqref{diagtwstko}. Since pullbacks of smooth 
stacks can be computed in prestacks, it suffices to produce the map $\iota$ objectwise, at the level of prestacks. 
Fix a smooth manifold $M$ and consider the triple of maps 
\begin{align}
& f_M:\Omega_{\rm cl}^1(M)  \xymatrix{\ar[r]&} {\rm Pic}_{\Omega^*\left(-;\pi_*(\op{KO})\right)}(M)\;, 
\label{drhmpk}
\\
& d_M:  \BB \RR^{\times}_{{\rm fl}\text{-}\nabla}(M) \xymatrix{\ar[r]&}
 {\rm Pic}_{H(\Omega^*(-;\pi_*(\op{KO})))}(M)\;, 
 \label{drmeqmp}
\\
&t_M:\BB\ZZ/2(M) \xymatrix{\ar[r]&} {\rm Pic}_{\delta(\op{KO})}(M)\;. 
\label{tpmapk}
\end{align}
The maps \eqref{drhmpk} and \eqref{tpmapk} are those described in diagram \eqref{diagtwstko} evaluated 
on $M$, namely ${\bf 2}$ and ${\bf 9}$, respectively. The map
 $d_M$ corresponds to the map induced by ${\bf 5}$ in \eqref{diagtwstko} which we now define explicitly. 
 By descent, an object in $\BB \RR^{\times}_{{\rm fl}\text{-}{\nabla}}(M)$ can be identified with a closed 1-form 
 $\mathcal{A}_{\alpha}$ defined on every open set of a fixed good open cover $\{U_{\alpha}\}$ of $M$. Then $d_M$ 
 assigns each $\mathcal{A}_{\alpha}$ to the sheaf of spectra $H(\Omega^*(-;\mathcal{L},d+\mathcal{A}_{\alpha}))$ 
 on $U_{\alpha}\subset M$. An edge in the stack corresponds to a gauge transformation 
 $g_{\alpha\beta}:\mathcal{A}_{\alpha}\mapsto \mathcal{A}_{\beta}+g^{-1}_{\alpha\beta}dg_{\alpha\beta}$ on 
 intersections and $d_M$ associates such an edge to the corresponding automorphism 
$$
\xymatrix{
g_{\alpha\beta}:H(\Omega^*(-;\mathcal{L},d+\mathcal{A}_{\alpha}))\ar[r]&
 H(\Omega^*(-;\mathcal{L},d+\mathcal{A}_{\beta}+g^{-1}_{\alpha\beta}dg_{\alpha\beta}))\;,
 }
 $$
which is induced by the map of complexes sending a form $\omega_{\alpha}\mapsto g_{\alpha\beta}\omega_{\beta}$. 
Descent implies that $d_M$ does not depend on the fixed cover $\{U_{\alpha}\}$ and this data is clearly natural with 
respect to pullback by smooth maps -- hence defines a morphism of smooth stacks.

We would like to produce the map $\iota:\BB(\ZZ/2)_{\nabla}\to \widehat{\rm Tw}_{\widehat{\op{KO}}}$ 
using the universal property of the pullback. Thus it suffices to send the homotopy in the pullback defining 
$\BB(\ZZ/2)_{\nabla}(M)$ to a homotopy in  ${\rm Pic}_{\delta(\op{KO})}(M)$. 
Again, utilizing descent, one finds that an edge in $\BB(\ZZ/2)_{\nabla}(M)$ can be identified with the following data 
(as in \cite{GS4}\cite{GS6}):
\begin{enumerate}[{\bf 1.}]
\vspace{-1mm} 
\item A choice of $\RR^{\times}$-valued smooth function $g_{\alpha}$ such that 
\(
 \label{ccocdat1} g_{\alpha}^{-1}dg_{\alpha}=\mathcal{A}_{\alpha}\;.
 \)
 
\vspace{-1mm} 
\item A choice of {\v C}ech cocycle $t_{\alpha\beta}\in \ZZ/2$, representing the topological twist $\sigma_1$, 
such that 
\(  \label{ccocdat2}g_{\alpha}g_{\beta}^{-1}=t_{\alpha\beta}\;.\)
\end{enumerate}
\vspace{-1mm} 
To these data, we must assign an edge in ${\rm Pic}_{H(\pi_*(\op{KO})\otimes \RR)}(M)$ connecting 
$H(\pi_*(\op{KO})_{\sigma_1}\otimes \RR)$ and the twisted de Rham complex 
$H(\Omega^*(-;\mathcal{L}),d+\mathcal{A})$. 
For this, we assign to $g_{\alpha}$ the local equivalence induced by the quasi-isomorphism of complexes
$$
g_{\alpha}\cdot:\pi_*(\op{KO})\otimes \RR\vert_{U_{\alpha}} 
\overset{\ker(d)}{\simeq} 
(\Omega^*(-;\pi_*(\op{KO})\otimes \RR),d)\vert_{U_{\alpha}}
\xymatrix{\ar[r]^{\simeq} &}
(\Omega^*(-;\mathcal{L}),d+\mathcal{A})\vert_{U_{\alpha}}
$$
where $g_{\alpha}=e^{f_{\alpha}}$ acts by the product $g_{\alpha}\cdot \omega_{\alpha}$. On intersections, we assign 
the {\v C}ech cocycle $t_{\alpha\beta}$ to the automorphism 
$$
t_{\alpha\beta}:\pi_*(\op{KO})\otimes \RR\vert_{U_{\alpha\beta}}\simeq 
\pi_*(\op{KO})\otimes \RR\vert_{U_{\alpha\beta}}\;,
$$
which is the transition data for the local system determined by $\sigma_1$. Passing to spectra via the Eilenberg-MacLane 
functor $H$, this is precisely the data necessary to specify an equivalence of sheaves of invertible module spectra
\(\label{eqshmsp}
d_{M}:H(\pi_*(\op{KO})\otimes \RR)_{\sigma_1}
\xymatrix{\ar[r]^{\simeq} &}
 (\Omega^*(-;\mathcal{L}),d+\mathcal{A})\;,
\)
or equivalently, an edge $d_{M}:\Delta[1]\to {\rm Pic}_{H(\pi_*(\op{KO})\otimes \RR)}(M)$ connecting the source 
and target of \eqref{eqshmsp}. Since the triple $(f_M,d_M,t_M)$ and the homotopy prescribed by the data 
\eqref{ccocdat1} and \eqref{ccocdat2} are natural with respect to pullback by smooth maps, it follows that we 
have produced the required homotopy and the universal property produces the map $\iota$. The rest of the claim follows 
immediately by construction.
\endofproof

The following example is of central importance and is a graded extension of one of the main 
examples in the construction of twisted Deligne cohomology in \cite{GS4}. Here we also see
that, even though the twist is torsion, differential refinement still captures interesting phenomena.

\begin{example}[Graded bundle of M\"obius strips]
Let $\mathcal{M}\to S^1$ be the "graded" M\"obius bundle $\mathcal{M}=\bigoplus_{i}M_{i}$ with $M_i$ a copy 
of the M\"obius bundle in degree $4i$. The trivial connection $\nabla=d$ has monodromy in $\ZZ/2\subset \RR^{\times}$. 
Consider the twist for $\op{KO}$ given by the classifying map $\sigma_1:S^1\to K(\ZZ/2,1)$ and let 
$\hat{\sigma}_1:S^1\to \BB(\ZZ/2)_{\nabla}$ be the differential refinement, defined via the universal property 
of the pullback via the following map and homotopy.

\begin{enumerate}
\item The map $0:S^1\to \Omega^1_{\rm cl}$ chooses the trivial connection on $\mathcal{M}$. 
\vspace{-0mm} 
\item The homotopy in $\BB \RR^{\times}_{{\rm fl}\text{-}\nabla}$ given by a choice of constant 
$\RR^{\times}$-valued functions $r_{\alpha}$, defined locally with respect to some cover $\{U_{\alpha}\}$, 
such that on intersections
$$
t_{\alpha\beta}=r_{\alpha}r_{\beta}^{-1}\in \ZZ/2\;,
$$
is a {\v C}ech cocycle representing the topological twist $\sigma_1$. 
\end{enumerate}

\vspace{-1mm}
\noindent Then the corresponding twist of $\widehat{\op{KO}}$ gives rise to a local system $\pi_*(\op{KO})_{\sigma_1}$, 
specified by the topological twist, the twisted de Rham complex $(\Omega^*(-;\mathcal{M}),d)$, and a twisted 
de Rham equivalence 
$$
d:H(\pi_*(\op{KO}_{\sigma_1})\otimes \RR)
\xymatrix{\ar[r]^{\simeq} &}
H((\Omega^*(-;\mathcal{M}),d))\;.
$$
In particular, $d$ induces a twisted de Rham isomorphism
$$H^*_{\rm dR}(S^1;\mathcal{M})\cong H^*(S^1;\pi_*(\op{KO})_{\sigma_1}\otimes \RR)\;.$$
\end{example}

\begin{remark}[Twists of degree 2]\label{Rem-twdeg2} 
In analogy with twists of degree 1 (see diagram \eqref{Bz2tw}), 
it is natural to think of refined twists classified by the smooth stack arising from the pullback 
$$
\xymatrix@C=3em@R=1.3em{
\BB^2(\ZZ/2)_{\nabla}\ar[r]\ar[d] & \Omega^2_{\rm cl}\ar[d]
\\
\BB^2\ZZ/2\ar[r] & \BB^2\RR^{\times}_{\nabla}\;.
}
$$
However, in this degree the rationalization map induces the trivial map on units and the differential form 
portion of  the twists do not serve to twist the corresponding de Rham complex (i.e. the rationalization does 
not have twists in degree 2).  Hence, there might seem that nothing to be gained by considering the geometry in this case. 
While an aspect of this is made precise by Proposition \ref{dftwtptw}, we will see later in Sec. \ref{Sec-AHSS-hatKO}
 that the geometry does get affected by torsion in a quite subtle way. 
\end{remark}

\begin{proposition}
[Space of determinantal differential twists for KO] \label{dftwtptw}
There is a map $\iota^{\prime}:\BB^2 \ZZ/2\into \widehat{\rm Tw}_{\widehat{{\rm KO}}}$, whose class in the 
homotopy category is a monomorphism, refining 
the inclusion of the factor $i^{\prime}:K(\ZZ/2,2)$ $\into B{\rm GL}_1({\rm KO})$. Moreover, for any
differential twist $\hat{\sigma}:M\to  \widehat{\rm Tw}_{\widehat{{\rm KO}}}$ refining a topological twist 
$\sigma:M\to K(\ZZ/2,2)\into B{\rm GL}_1(\op{KO})$, we have (at the level of homotopy classes)
$$
[\hat{\sigma}]-[i'(\sigma)]=k([\overline{\sigma}])\;,
$$ 
where $k:{\rm fib}(\mathcal{R})\to \widehat{{\rm Tw}}_{\widehat{\op{KO}}}$ is the canonical map induced 
by taking the fiber of the top horizontal map in \eqref{Tw-KO} and $[\overline{\sigma}]$ represents a class in 
$\pi_0\map(M,{\rm fib}(\mathcal{R}))$. In particular, the corresponding twisted de Rham complex is 
equivalent to the untwisted complex.

\end{proposition}
\theproof
Consider the map 
$$
\xymatrix{
j:\BB^2\ZZ/2\ar[r]& \delta(\vert \BB^2 \ZZ/2 \vert)\simeq \delta(B^2\ZZ/2)
\; \ar@{^{(}->}[r]& \delta(B{\rm GL}_1(\op{KO})) \; \ar@{^{(}->}[r]& {\rm Pic}_{\delta(\op{KO})}},
$$
where $\vert \cdot \vert$ denotes the geometric realization of the stack and $\delta$ denotes the constant 
stack functor. The first map in the composition $j$ is the unit of the adjunction 
$\eta:\mathbb{1}\to \delta\circ \vert \cdot \vert$. 
Post-composition of $j$ with the rationalization 
$\wedge H\RR:{\rm Pic}_{\delta(\op{KO})}\to {\rm Pic}_{H(\pi_*(\op{KO}))\wedge H\RR)}$
 is a null-homotopic map. By the universal property, a choice of such a nullhomotopy gives rise to an induced map 
$i:\BB^2\ZZ/2\to \widehat{{\rm Tw}}_{\widehat{\op{KO}}}$,
whose post-composition with the canonical map 
$\widehat{{\rm Tw}}_{\widehat{\op{KO}}}\to {\rm Pic}_{\Omega^*(-;\pi_*(\op{KO}))}$ 
is null-homotopic and whose post-composition with 
$\mathcal{I}:\widehat{\rm Tw}_{\widehat{{\rm KO}}}\to {\rm Pic}_{\delta({\rm KO})}$ 
is equivalent to $i$. Now the homotopy functor
$$
{\rm Ho}(\delta):{\rm Ho}(\infty\text{-}\mathcal{G}{\rm pd})
\xymatrix{\ar[r]&}
 {\rm Ho}(\sh_{\infty}({\mathcal{C}\mathrm{art}\mathcal{S}\mathrm{p}}))
$$
preserves monomorphisms, since it is a right adjoint, and 
$$
\mathcal{I}i:\delta(K(\ZZ/2,1))\longrightarrow 
\delta({\rm GL}_1(\op{KO}))\longrightarrow {\rm Pic}_{\delta(\op{KO})}
$$ 
is a monomorphism in the homotopy category (given by inclusion of the factor). Hence so is $i$. 

Now let $\sigma:M\to \BB^2\ZZ/2\simeq \delta(K(\ZZ/2,2))\into \delta(B{\rm GL}_1(\op{KO}))\into {\rm Pic}_{\delta(\op{KO})}$ 
be a twist for $\op{KO}$ and let $\hat{\sigma}:M\to \widehat{{\rm Tw}}_{\widehat{{\op{KO}}}}$ be any refinement of $\sigma$. 
By definition, it follows that $\mathcal{I}([\hat{\sigma}]-[\hat{i}(\sigma)])=\mathcal{I}([\hat{\sigma}])-\mathcal{I}([\hat{\sigma}])=0$. 
The long exact sequence on homotopy associated to the fiber sequence 
$$
\xymatrix{
\map(M,{\rm fib}(\mathcal{R}))\ar[r]^-{k} & 
\map(M,\widehat{{\rm Tw}}_{\widehat{{\op{KO}}}})\ar[r]^-{\mathcal{I}} &
\map(M,{\rm Pic}_{\delta(\op{KO})})
},
$$
in particular for $\pi_0$, gives the existence of $\bar{\sigma}$. 
\endofproof

Combining Propositions \ref{dftwtptw} and \ref{sttwdeg1} gives the following main result of this section.

\begin{theorem}[Stack of twists of $\widehat{\rm KO}$-theory]
\label{mainthtw}
There is a morphism of smooth stacks
$$
\xymatrix{
\iota:\BB(\ZZ/2)_{\nabla}\times \BB^2\ZZ/2\ar[r] & \widehat{{\rm Tw}}_{\widehat{\op{KO}}}
}
$$
refining the inclusion $i:K(\ZZ/2,1)\times K(\ZZ/2,2)\into B{\rm GL}_1(\op{KO})$, whose restriction to the factor 
$\BB^2\ZZ/2$ is  a monomorphism in the homotopy category. Moreover, any twist 
$\hat{\sigma}:M\to \widehat{\rm Tw}_{\widehat{\op{KO}}}$ 
which factors through $\iota$ corresponds to a twisted differential cohomology theory which fits into the pullback
 \footnote{This pullback is taken in the slice of the tangent topos of parametrized sheaves of spectra over $M$.}
$$
\xymatrix@R=1.3em@C=3em{
\widehat{\op{KO}}_{\hat{\sigma}}\ar[r]\ar[d] & H(\tau_{\geq 0}\Omega^*(-;\mathcal{L}))\ar[d]
\\
\delta(\op{KO}_{\sigma})\ar[r] & H(\Omega^*(-;\mathcal{L}))\;,
}
$$
where the homotopy can be identified with an equivalence 
$$
d:H(\pi_*(\op{KO})_{\sigma_1}\otimes \RR)
\xymatrix{\ar[r]^{\simeq} &}
H(\Omega^*(-;\mathcal{L}))\;,
$$
the complex $\Omega^*(-;\mathcal{L})$ is defined as in \eqref{twfrmdg1} and 
$\sigma=(\sigma_1,\sigma_2):M\to K(\ZZ/2,1)\times K(\ZZ/2,2)$ is the underlying topological twist.
\end{theorem}

\subsection{Properties of twisted differential KO-theory}
\label{Sec-prop}

Having discussed the twists of $\widehat{\op{KO}}$, we now turn our attention to the basic properties 
of the twisted theory.  

\begin{proposition}
[Properties of twisted differential KO-theory]\label{prop-prop}
Let $\hat{\sigma}:M\to \BB(\ZZ/2)_{\nabla}\times \BB^2\ZZ/2$ be a twist for $\widehat{{\rm KO}}$, as in 
Theorem \ref{mainthtw}. Then the twisted differential cohomology theory $\widehat{{\rm KO}}_{\sigma}$ 
satisfies the following properties:

\item {\bf (i)}  \emph{(Integration and curvature)}
We have maps
\( \label{IandI}
\mathcal{I}_{\sigma}:\widehat{{\rm KO}}^*_{\sigma}(M)\longrightarrow {\rm KO}^*_{\sigma}(M)
\qquad
\text{and}
\qquad 
 \mathcal{R}_{\sigma}:\widehat{{\rm KO}}^*_{\sigma}(M)\longrightarrow \Omega^*(M;\mathcal{L})\;,
\)
which are natural in the sense that, for any smooth map $f:N\to M$, we have commutative diagrams 
\(\label{cmdiagtwsp}
\hspace{-3mm}
\xymatrix{
\widehat{{\rm KO}}^*_{\sigma}(M)\ar[rr]\ar[d]_-{\mathcal{I}_{\sigma}}  &&
 \widehat{{\rm KO}}^*_{f^*\sigma}(N)\ar[d]^{\mathcal{I}_{f^*\sigma}} 
\\
{\rm KO}^*_{\sigma}(M)\ar[rr] && \widehat{{\rm KO}}^*_{f^*\sigma}(N)
}
\qquad 
\text{and}
\qquad
\xymatrix@R=2.7em{
\widehat{{\rm KO}}_{\sigma}(M)\ar[rr]\ar[d]_-{\mathcal{R}_{\sigma}}  && 
\widehat{{\rm KO}}_{f^*\sigma}(N)\ar[d]^{\mathcal{R}_{f^*\sigma}} 
\\
\Omega^*(M;\mathcal{L})\ar[rr] && \Omega^*(N;\mathcal{L})
\;.}
\)
\item{\bf (ii)} \emph{(Secondary forms)} There is a canonical homomorphism of \emph{groups}
$$
a_{\sigma}:\Omega^{*-1}(M;\mathcal{L})/{\rm Im}(d)\longrightarrow
 \widehat{{\rm KO}}^*_{\sigma}(M)
$$
which satisfies
$$a_{\sigma}(\omega)\cdot x= a_{\sigma}(\omega\wedge \mathcal{R}_{\sigma}(x))
\qquad 
\text{and}
\qquad
 \mathcal{R}_{\sigma}\circ a_{\sigma}=d\;.
 $$
\item {\bf (iii)} \emph{(Diamond)}  We have the twisted differential cohomology diamond diagram 
\(
\hspace{-.6cm}
\label{twkodfdiam}
\xymatrix @C=-5pt @!C{
 &\Omega^{*-1}(M;\mathcal{L})/{\rm im}(d) \ar[rd]_{a_{\sigma}}\ar[rr]^{d} 
 & &  \Omega_{\rm fl}^*(M;\mathcal{L})\ar[rd] &  
\\
H_{\rm dR}^{*-1}(M;\mathcal{L})\ar[ru]\ar[rd] & & 
{\widehat{{\rm KO}}_{\sigma}^*(M)}\ar[rd]^{\mathcal{I}_{\hat{\sigma}}} \ar[ru]_{\mathcal{R}_{\sigma}}& &  
H_{\rm dR}^*(M;\mathcal{L}) \;.
\\
&{\rm KO}_{\sigma}^{*-1}(M;U(1))\ar[ru]^-j\ar[rr]^{\beta} & & {\rm KO}_{\sigma}^*(M) \ar[ru] &
}
\) 
where $\mathcal{L}$ is the graded line bundle with connection, classified by the degree 1 twist
$$\hat{\sigma}_1:M\to \BB(\ZZ/2)_{\nabla}\;.$$
\end{proposition}
\theproof
{\bf (i)} From Theorem \ref{mainthtw}, the twists of the form
$$
\hat{\sigma}:M\longrightarrow
 \BB(\ZZ/2)_{\nabla}\times \BB^2\ZZ/2\longrightarrow
  \widehat{\rm Tw}_{\widehat{{\rm KO}}}\;,
$$ 
correspond to the sheaves of spectra defined by the pullback diagrams
\(\label{shsptw}
\xymatrix@C=4em{
{\rm diff}({\rm KO}_{\sigma},d,\mathcal{L})\ar[r]^-{\mathcal{R}_{\sigma}}\ar[d]^-{\mathcal{I}_{\sigma}} 
& H(\tau_{\geq 0}\Omega^*(-;\mathcal{L}))\ar[d]
\\
{\rm KO}_{\sigma} \ar[r] & H(\Omega^*(-;\mathcal{L}))
}
\)
of sheaves of spectra on $M$. We take $\mathcal{I}_{\sigma}$ and $\mathcal{R}_{\sigma}$ in \eqref{IandI} 
to be induced by the corresponding maps in \eqref{shsptw}. Considering the shifts of the spectrum $\op{KO}$, we get 
corresponding diagrams in all degrees. The commutative diagrams \eqref{cmdiagtwsp} follow immediately by 
functoriality of the above construction. 

\noindent {\bf (ii)} The subcategory 
$T_{M}(\sh_{\infty}({\mathcal{C}\mathrm{art}\mathcal{S}\mathrm{p}}))\into 
T(\sh_{\infty}({\mathcal{C}\mathrm{art}\mathcal{S}\mathrm{p}}))$ 
of the tangent $\infty$-topos, 
given by the full sub $\infty$-category of spectrum objects over $M$ is a stable $\infty$-category. \footnote{Indeed, via the $\infty$-Grothendieck construction, this can be identified with the stabilization of the functor category 
${\rm Fun}(X,\sh_{\infty}({\mathcal{C}\mathrm{art}\mathcal{S}\mathrm{p}}))$}
 The objects in this category are identified with sheaves of spectra on $M$.  In any stable $\infty$-category a square is Cartesian if 
and only if the induced map on fibers over the zero object is an equivalence (see \cite[Remark 7.1.12]{Hov}). 
Using standard techniques in homological algebra, one finds that we have a fiber sequence of sheaves of spectra on $M$
\(\label{fibseqsh}
\xymatrix{
H(\tau_{<0}\Omega^{*}(-;\mathcal{L}))\ar[r] &
H(\tau_{\geq 0}\Omega^{*}(-;\mathcal{L}))\ar[r] &
 H(\Omega^{*}(-;\mathcal{L}))
 }.
\)
Then, by stability, it follows from the sequence \eqref{fibseqsh} and the pullback diagram \eqref{shsptw}, that we have a 
fiber sequence 
$$
\xymatrix{
 H(\tau_{<0}\Omega(-;\mathcal{L}))\ar[r] & {\rm diff}({\rm KO}_{\sigma},d,\mathcal{L})\ar[r] & {\rm KO}_{\sigma}
 }
 $$
in $T_{M}(\sh_{\infty}({\mathcal{C}\mathrm{art}\mathcal{S}\mathrm{p}}))$. 
By the Poincar\'e Lemma (which holds for forms with values in $\mathcal{L}$), 
the global sections of the latter are identified with 
$\Omega^{-1}(M;\mathcal{L})/{\rm im}(d)$. We define $a_{\sigma}$ as the map induced by the left vertical 
map in the above diagram. Considering the shifts of the spectrum $\op{KO}_{\sigma}$ gives the construction in all 
degrees. The formula relating the products follows from \cite[Proposition 5.1]{BN}. 

\noindent {\bf (iii)}  The existence of the diamond or hexagon follows from the construction in \cite{BNV}, which holds 
in any stable $\infty$-category.  We apply the construction in the subcategory 
$T_{M}(\sh_{\infty}({\mathcal{C}\mathrm{art}\mathcal{S}\mathrm{p}}))$ to the 
sheaf of spectra ${\rm diff}(\op{KO}_{\sigma},d,\Omega^*(-;\mathcal{L}))$. The only possibly nontrivial identification 
is the spectrum representing the twisted flat theory $\op{KO}_{\sigma}^{*-1}(-;U(1))$. Considering the fiber sequence 
$
\mathbb{S}\to \mathbb{S}\RR\to \mathbb{S}\RR/\ZZ
$,
we have an associated fiber sequence of locally constant sheaves of spectra over $M$
\(\label{fibstwsp}
\xymatrix{
\delta(\op{KO}_{\sigma})\ar[r] &
 \delta(\op{KO}_{\sigma})\wedge \mathbb{S}\RR\ar[r] &
\delta(\op{KO}_{\sigma})\wedge \mathbb{S}\RR/\ZZ
}.
\)
Since the twist is torsion, as in the proof of Proposition \ref{dftwtptw}, we see that the Pontrjagin character gives 
an identification
$$
\xymatrix{
\delta({\rm Ph}):\delta(\op{KO}_{\sigma})\wedge \mathbb{S}\RR\ar[r]^-{\simeq} &
 \delta(H(\pi_*(\op{KO})\otimes \mathcal{L}))
 }.
$$
The long exact sequence associated to \eqref{fibstwsp} gives rise to the fiber sequence
$$
\xymatrix{
\Sigma^{-1}_M\delta(\op{KO}_{\sigma}\wedge \mathbb{S}\RR/\ZZ)\ar[r] &
\delta(\op{KO}_{\sigma})\ar[r] &
 \delta(H(\pi_*(\op{KO})\otimes \mathcal{L}))
},
$$
which identifies the flat theory on the hexagon. The same argument as in {\bf (ii)}
yields the the lower portion of the  diamond diagram. The identification of the left square 
of the diamond proceeds by the dual argument (i.e., pushout vs. pullback).
\endofproof

In \cite{GS5} we constructed an Atiyah-Hirzebruch spectral sequence which works for any twisted differential 
cohomology theory. The filtration which gives rise to the spectral sequence is given by restricting a twist 
$\sigma:M\to \widehat{{\rm Tw}}_{\widehat{\mathcal{R}}}$ to elements of a good open cover $\{U_i\to M\}$ 
and considering the filtration of the {\v C}ech nerve of this cover by skeleta. For $\widehat{\op{KO}}$, the 
construction immediately gives the following.
\begin{proposition}[Twisted differential AHSS]
\label{spseqdfko}
Let $M$ be a compact smooth manifold and let 
$\hat{\sigma}=(\hat{\sigma_1},\sigma_2):M\to \BB(\ZZ/2)_{\nabla}\times \BB^2 \ZZ/2$ be a twist. 

\medskip
\noindent {\bf (i)} There is a half-plane spectral sequence with abutment 
$\widehat{\op{KO}}_{\hat{\sigma}}(M)$ with $E_2$-page given by 
$$
E_2^{p,q}=\left\{ \begin{array}{ccc}
\Omega_{{\rm fl}}^{0}(M; \mathcal{L}) && p=q=0,
\\
H^p(M;\mathcal{L}^{\delta}/\mathcal{Z}) && p>0, q<0, q\equiv -3 \!\!\mod 4
\\
H^p(M;\ZZ/2) && p>0, q<0, q\equiv -1,-2\!\!\mod 4
\\
0 && \text{otherwise}
\end{array}
\right.
$$ 
Note that, in bidegree $(0,0)$, we have twisted closed differential forms with coefficients in the flat graded 
line bundle $\mathcal{L}\to M$, with fiber $\pi_*({\rm KO})\otimes \RR \simeq \RR[\alpha,\alpha^{-1}]$, 
$\vert \alpha\vert=4$, determined by the factor $\hat{\sigma}_1$, The discrete bundle $\mathcal{L}^{\delta}$ 
is the local system associated to the flat connection and $\mathcal{Z}\into \mathcal{L}^{\delta}$ is the 
$\ZZ$-subbundle classified by topological twist $\sigma_1$. 

\medskip
\noindent {\bf (ii)} In the case where $\sigma_1=0$, the $E_2$-page takes the form 
$$
E_2^{p,q}=\left\{ \begin{array}{ccc}
\Omega_{{\rm cl}}^{0}(M; \RR[\alpha,\alpha^{-1}])  && p=q=0,
\\
H^p(M;\RR/\ZZ) && p>0, q<0, q\equiv -3 \!\!\mod 4
\\
H^p(M;\ZZ/2) && p>0, q<0, q\equiv -1,-2\!\!\mod 4
\\
0 && \text{otherwise}
\end{array}
\right.
$$ 
\end{proposition}

\subsection{The twisted differential Pontrjagin character} 
\label{Sec-tdPh}

In this section we define the twisted differential Pontrjagin character -- the differential refinement of the 
corresponding character in twisted KO-theory in Sec. \ref{Sec-tPh}. By \cite[Proposition 36]{GS-KO}, 
there is a differential refinement of the Pontrjagin character, which is induced by a morphism of sheaves of spectra
\(\label{dfphmp}
\widehat{\rm Ph}:{\rm diff}({\rm KO},{\rm Ph},\pi_*({\rm KO}))
\xymatrix{\ar[r]&}
 {\rm diff}(H(\pi_*({\rm KO})\otimes \QQ),i,\pi_*({\rm KO})\otimes \RR)\;.
\)
Moreover, the induced natural transformation 
$\widehat{{\rm Ph}}:\widehat{\op{KO}}\to \widehat{H}(-;\pi_*(\op{KO})\otimes \QQ)$ 
is the unique transformation which refines the topological Pontrjagin character and the Pontrjagin character 
form \cite[Corollary 37]{GS-KO}. The morphism \eqref{dfphmp} induces a morphism on the corresponding 
stack of twists as follows. For each invertible module $\mathcal{L}$ over $\op{KO}$, form the module 
$\mathcal{L}\wedge_{\op{KO}}H(\pi_*(\op{KO})\otimes \RR)$, where $H(\pi_*(\op{KO})\otimes \RR)$ 
is viewed as a $\op{KO}$ -module via the morphism of ring spectra ${\rm Ph}$. This is an invertible module over 
$H(\pi_*(\op{KO})\otimes \RR)$, with inverse $\mathcal{L}^{-1}\wedge_{\op{KO}}H(\pi_*(\op{KO})\otimes \RR)$,
 and the construction is functorial. We then have an obvious commutative diagram
$$
\xymatrix@C=4em{
{\rm Pic}_{\Omega^*(-;\pi_*(\op{KO}))}\ar@{=}[rr] \ar[d]&& {\rm Pic}_{\Omega^*(-;\pi_*(\op{KO}))}\ar[d]
\\
{\rm Pic}_{H(\Omega^*(-;\pi_*(\op{KO}))}\ar@{=}[rr] && {\rm Pic}_{H(\Omega^*(-;\pi_*(\op{KO}))}
\\
{\rm Pic}_{\delta(\op{KO})}\ar[u]\ar[rr]^-{\delta(\wedge_{\op{KO}}H(\pi_*(\op{KO})\otimes \QQ))} && 
{\rm Pic}_{\delta(H(\pi_*(\op{KO})\otimes \RR))}\;. \ar[u]
}
$$
From the universal property of the pullback, we thus have the following.

\begin{proposition}[The differential Pontrjagin character] \label{prop-diffpontch}
The differential Pontrjagin character induces a morphism of smooth stacks
$$
\widetilde{{\rm Ph}}:\widehat{{\rm Tw}}_{\widehat{\op{KO}}}\longrightarrow 
\widehat{{\rm Tw}}_{\widehat{H}(-;\pi_*(\op{KO})\otimes \QQ)}\;.
$$
\end{proposition} 

\medskip
We now turn our attention to the twisted Pontrjagin character. We use a similar method as that in Sec.
\ref{Sec-tPh} to prove the existence and uniqueness. Just as twisted topological theories are represented by 
objects in the tangent $\infty$-topos modeled on $\infty$-groupoids, twisted differential theories are represented 
by objects in the tangent $\infty$-topos modeled on smooth $\infty$-groupoids 
$\sh_{\infty}(\cartsp)$ 
(see \cite{GS5} for detailed discussion). The axiom of descent holds for $T(\sh_{\infty}(\cartsp))$ and we can 
apply the same local-to-global technique, as we did in Sec. \ref{Sec-tPh}. 

\begin{theorem}[The twisted differential Pontrjagin character]\label{thm-twdifPo} 
Let $M$ be a smooth manifold and let $\hat{\sigma}=(\hat{\sigma}_1,\sigma_2):M\to \BB(\ZZ/2)_\nabla\times \BB^2\ZZ/2$ be 
a twist for $\widehat{\op{KO}}$ with underlying topological twist $\sigma=(\sigma_1,\sigma_2)$. Let $\mathcal{L}\to M$ 
be the real graded line bundle with connection defined by the twist $\hat{\sigma}_1$. There is a morphism of corresponding sheaves 
of spectra on $M$
$$
\xymatrix{
\widehat{{\rm Ph}}_{\sigma}:{\rm diff}({\rm KO}_{\sigma},d,\mathcal{L}) \ar[r]& 
  {\rm diff}(H(\pi_*({\rm KO})_{\sigma_1}\otimes \QQ),d^{\prime},\mathcal{L})
  },
  $$
called the \emph{twisted differential Pontrjagin character}, which refines the underlying topological Pontrjagin character. 
Here we regard $\hat{\sigma}$ as a twist for $H(\pi_*({\rm KO})\otimes \QQ)$ via rationalization as described 
in Sec. \ref{Sec-tPh}.
\end{theorem}

\theproof
Fix a good open cover $\{U_{\alpha}\}$ of $M$. Pulling back the sheaf of spectra
$\widehat{\op{KO}}_{\sigma}\to M$ by the inclusion $i_{\alpha}:U_{\alpha}\into M$, gives a sheaf of spectra 
$\widehat{\op{KO}}_{i^*_{\alpha}\sigma}$ over $U_{\alpha}$. Since $U_{\alpha}$ is contractible, 
Proposition 4 in \cite{GS5}
implies that we have an equivalence 
$$
f_{\alpha}: {\rm diff}\big({\rm KO}_{i^*_{\alpha}\sigma},i^*_{\alpha}d,i^*_{\alpha}\mathcal{L}\big)
\simeq 
{\rm diff}({\rm KO},{\rm Ph},\pi_*({\rm KO}))
$$
of sheaves of spectra over $U_{\alpha}$. Similarly, we have local trivializations 
$$
g_{\alpha}:{\rm diff}\big(H(\pi_*({\rm KO})\otimes \QQ)_{i^*_{\alpha}\sigma},
i^*_{\alpha}d^{\prime},i^*_{\alpha}\mathcal{L}\big)
\simeq 
{\rm diff}\big(H(\pi_*({\rm KO})\otimes \QQ),d^{\prime},\pi_*({\rm KO})\otimes \RR\big)
$$
over each $U_{\alpha}$. We can locally define the twisted differential Pontrjagin character,
${\widehat{{\rm Ph}}_{i^*_{\alpha}\sigma}}$,
as the clockwise composition
$$
\hspace{-2mm}
\xymatrix@C=2.5em{
{\rm diff}({\rm KO}_{i^*_{\alpha}\sigma},i^*_{\alpha}d,i^*_{\alpha}\mathcal{L})\ar[d]_{\widehat{{\rm Ph}}_{i^*_{\alpha}\sigma}}
\ar[rr]_-{\simeq}^{f_{\alpha}} && 
{\rm diff}({\rm KO},{\rm Ph},\pi_*({\rm KO}))\ar[d]^-{\widehat{{\rm Ph}}}
\\
{\rm diff}(H(\pi_*({\rm KO})\otimes \QQ)_{i^*_{\alpha}\sigma^*{\rm Ph}},i^*_{\alpha}d^{\prime},i^*_{\alpha}\mathcal{L})&& 
\ar[ll]^-{\simeq}_-{g^{-1}_{\alpha}} {\rm diff}(H(\pi_*({\rm KO})\otimes \QQ),i,\pi_*({\rm KO})\otimes \RR).
}
$$
Since the differential Pontrjagin character induces a morphism of smooth stacks (see Prop. \ref{prop-diffpontch}), 
it follows that the local twisted Pontrjagin character maps $\widehat{{\rm Ph}}_{i^*_{\alpha}\sigma}$ are compatible 
 with the transition data for the bundle of spectra. By descent it follows that the local morphisms 
 $\widehat{{\rm Ph}}_{i_{\alpha}^*\sigma}$ induce a global morphism which is unique for this data, up to a contractible choice. 
\endofproof

\section{The Atiyah-Hirzebruch spectral sequence for twisted differential KO-theory}
\label{Sec-AHSS}

In this section we develop our computational tools, in the form of the AHSS. It turns out that
we need to devote a lot of attention to the twisted topological theory before being able to 
discuss the twisted differential theory.

\subsection{AHSS for twisted KO-theory}
\label{Sec-tAHSS} 

The goal of this section is to identify some of the differentials in the spectral sequence for twisted $\op{KO}$-theory. 
Interestingly, it appears that these have not yet been calculated, even for the $E_2$-page of the spectral sequence. 
Because we are also considering twists of degree 1, the $E_2$-page of the AHSS will contain cohomology with local 
coefficients. Thus we have a real analogue of the statements in \cite{AS2} for the complex case, and extending the real 
case from having just a degree two twist in \cite{Ros}. 

\medskip
Given $\sigma=(\sigma_1,\sigma_2):X\to K(\ZZ/2,1)\times K(\ZZ/2,2)$, the AHSS for twisted $\op{KO}$-theory 
takes the following form
\(
E_2^{p,q}=H^p(X;\pi_{q}(\op{KO})_{\sigma_1})\Longrightarrow \op{KO}^{p+q}_\sigma(X)\;,
\)
where $\pi_{q}(\op{KO})_{\sigma_1}$ is the local coefficient system specified by the degree 1 
twist $\sigma_1:X\to K(\ZZ/2,1)$ and the  canonical action of the units $\ZZ/2$ on the graded ring 
$\pi_{*}(\op{KO})$. Note that this action is trivial for $\pi_{q}(\op{KO})\cong \ZZ/2$  with $q=1,2\mod 8$. 


\medskip
Consider the canonical bundle map $\pi:K(\ZZ/2,p)\sslash \ZZ/2\longrightarrow B\ZZ/2$, where the quotient denotes 
the homotopy orbits of the canonical $\ZZ/2$-action. By a standard representability arguments, the homotopy classes of maps 
$K(\ZZ/2,p)\sslash \ZZ/2\longrightarrow K(\ZZ/2,q)\sslash\ZZ/2$ covering $\pi$ are in bijective correspondence with cohomology 
operations with local coefficient system specified by a map $\sigma_1:X\to B\ZZ/2$. Denote the set of such homotopy classes by 
$H^q(K(\ZZ/2,p)\sslash\ZZ/2; \mathcal{Z}_2)$, where $\mathcal{Z}_2$ indicates the local coefficient system. According to the 
main theorem of  Gitler in \cite{Git}, the quotient map $i:K(\ZZ/2, p)\longrightarrow K(\ZZ/2, p)\sslash\ZZ/2$ 
induces an isomorphism of rings
$$
\xymatrix{
i^*:H^*(K(\ZZ/2,p)\sslash\ZZ/2; \mathcal{Z}_2)\ar[r]^-{\cong}& H^*(K(\ZZ/2,p);\ZZ/2)
}.
$$
Hence, taking local coefficients has no effect on the possible natural operations. 
We immediately get the following. 

\begin{lemma}[Operations with local coefficients]
The natural operations in cohomology with local coefficients for spaces equipped with maps 
$\sigma:X\to K(\ZZ/2,1)\times K(\ZZ/2,2)$ raising degree by 2 are thus in 
bijective  correspondence with  the cohomology group
$$
H^{p+2}\big(K(\ZZ/2,1)\times K(\ZZ/2,2)\times K(\ZZ/2,p);\ZZ/2\big)\;.
$$
\end{lemma}
From the K\"unneth formula, one finds that the only possibility is combinations of $Sq^2$, $Sq^1$, $w_1$ and $w_2$, the 
latter corresponding to the generators of $H^1(K(\ZZ/2,1);\ZZ/2)$ and $H^2(K(\ZZ/2,2);\ZZ/2)$. Similar considerations 
also hold for operations from $\ZZ$-coefficients to $\ZZ/2$-coefficients and $\ZZ/2$-coefficients to $\ZZ$-coefficients. 
Since the first non-vanishing differentials are necessarily stable cohomology operations, it follows from these 
considerations that we have the following. 

\begin{remark}[Possible form of the differentials] \label{rem-pos}  
The differentials take the general form
\begin{align}
d^{p,-8t}_2 &=aSq^2\circ r+b\sigma_1\cup Sq^1\circ r+c \sigma_1^2\cup r+d\sigma_2\cup r\;, 
\label{dfeqns1}
\\
d^{p,-8t-1}_2 &=a^{\prime}Sq^2+b^{\prime}\sigma_1\cup Sq^1+c^{\prime} \sigma_1^2\cup +d^{\prime} \sigma_2\cup\;,
\label{dfeqns2}
\\
d^{p,-8t-2}_2 &=a^{\prime\prime}\beta\circ Sq^2+d^{\prime\prime}\beta\circ \sigma_2\cup\;,
 \label{dfeqns3}
\end{align}
where $r:H^*(X;\ZZ)\to H^*(X;\ZZ/2)$ is the mod 2 reduction, $\beta:H^*(X;\ZZ/2)\to H^{*+1}(X;\ZZ)$ is the Bockstein 
homomorphism, and the coefficients are in $\ZZ/2$, and where the constants $a$, $a^{\prime}$, $a^{\prime\prime}$,
$b$, $b^{\prime}$,  $c$, $c^{\prime}$, $d$, $d^{\prime}$, and $d^{\prime\prime} \in \Z/2$ are to be determined.
\end{remark} 

Our work is thus reduced to identifying the coefficients. Note that the first three 
coefficients are $a=a^{\prime}=a^{\prime\prime}=1$, as these correspond to the untwisted case 
 \cite{Fu1}\cite{ABP}\cite{GS-KO}. 

\medskip
To identify the remaining coefficients, we will provide examples of spaces which isolate each term of the differentials. 
In each case, we either show that the differential cannot vanish, and thus the corresponding coefficient is 1, or that 
having a non-vanishing differential leads to a contradiction. 

\medskip
The real projective spaces  $\RR P^n$ and higher dimensional Klein bottles $K_n$
will be essential for our identification, but in order to use these spaces, we need to make some preliminary calculations. 
First, we use the main results of \cite{Fu1}
(see \cite{Aj} for a simplified treatment), 
which computes the $\widetilde{\op{KO}}$-theory of real projective spaces -- these are tabulated in the following 
table 
and we then compute the $\op{KO}$ groups and twisted $\op{KO}$ groups of Klein bottles.

\begin{table}[h]
\hspace{0cm}
{\small
\begin{tabular}{|c||c|c|c|c|c|c|c|c|}
\hline 
\diagbox{$i$}{$n$} &  $8r$ &  $8r+1$ & $8r+2$ &
$ 8r+3$ &
$8r+4$ & $8r+5$ & $8r+6$ & $8r+7$
\\ \hline \hline
$ 0$ & $(2^{4r})$ & $(2^{4r+1})$  & $(2^{4r+2})$  & $(2^{4r+2})$  & $(2^{4r+3})$ & $(2^{4r+3})$  & $(2^{4r+3})$  & $(2^{4r+3})$ 
\\ \hline
$1$  & $(2)$ & $(2)$ & $(2)$ & $(\infty)+(2)$  & $(2)$  & $(2)$ &$(2)$ & $(\infty)+(2)$
\\ \hline
$ 2$  & \!$(2)+(2)$\!  & $(2)$ & $(2)$& $(2)$&  $(2)$ & $(2)$ &  $(2)+(2)$ & \!\!$(2)+(2)+(2)$ \!\!\!
\\ \hline
$3$  & $(2)$ & $(\infty)$ & $0$ & $0$ & $0$ & $(\infty)$ & $(2)$ & $(2)+(2)$
\\ \hline
$4$ & $(2^{4r})$ & $(2^{4r})$  & $(2^{4r})$  & $(2^{4r})$  & $(2^{4r+1})$ & $(2^{4r+2})$  & $(2^{4r+3})$  & $(2^{4r+3})$ 
\\ \hline
$5$  & $0$ & $0$&  $0$ & $(\infty)$ & $0$ & $0$ & $0$&  $(\infty)$
\\ \hline
$6$  & $0$ & $0$&  $(2)$ & $(2)+(2)$ & $(2)$ &  $0$ & $0$&  $0$
\\ \hline
$7$  & $0$ &  $(\infty)$  &  $(2)$ & $(2)+(2)$  &  $(2)$ & $(\infty)$ & $0$&  $ 0$
\\ \hline
\end{tabular}
}
\label{korealprfu}
\caption{The groups $\widetilde{\op{KO}}^{-i}(\RR P^n)$. Here, $(k)$ denotes the cyclic group of order $k$.}
\end{table}

\begin{remark}[Utility of the twisted Thom isomorphism and real projective spaces]
\label{twkorpntiso}
One useful trick is to use the twisted Thom isomorphism in KO-theory, established in \cite[Prop. 4.1.7]{HJ},
in order to compute the twisted $\op{KO}$-theory of real projective spaces. In particular, if we take the twist 
$(w_1,w_2): \RR P^n\to K(\ZZ/2,1)\times K(\ZZ/2)$ corresponding to $w_1$ and $w_2$ of the canonical 
line bundle $\xi\to \RR P^n$, then $w_2$ vanishes since its degree exceeds the rank. From the identification 
${\rm Th}(\xi)\simeq \RR P^{n+1}$, we see that the twisted Thom isomorphism gives an isomorphism
$$
\widetilde{\op{KO}}_{w_1}^{-i}(\RR P^n) \cong \widetilde{\op{KO}}^{-i+1}(\RR P^{n+1})\;.
$$
Thus, we can read off the twisted $\op{KO}$-groups of real projective spaces from the above table, simply by 
referring to the  entry diagonally up and right from the entry corresponding to the untwisted group.
\end{remark}

\begin{remark}[Untwisting effect] \label{rem-untw}
Note that this untwisting effect has been observed in 
twisted K-theory in \cite{Kh} and in twisted Morava K-theory
in \cite{SY}. Thus classical calculations of the latter untwisted groups
$KO(TE)$, e.g. in \cite{ABP}\cite{Dob}, where $TE$ is the Thom space of a bundle $E$,
 can be recast as being twisted KO-theory calculations. 
 We will discuss these matters in detail elsewhere. 
 \end{remark} 

\medskip
We now turn to the $\op{KO}$-theory of Klein bottles. In \cite{Dav}, the cohomology groups of higher-dimensional 
Klein bottles were identified.  It is also tacitly assumed there that a combination of the AHSS and Mayer-Vietoris 
sequence could be used to compute the $\op{K}$-theory  and $\op{KO}$-theory of Klein bottles. While it seems that in 
principle one could do this, the calculations involved appear to us to be quite tedious. Therefore, we have chosen 
a different approach which involves an interesting interplay between twisted and untwisted cohomology. We also 
use the method to compute the twisted KO-theory.

\begin{proposition}[Twisted KO-theory of Klein bottles]
Consider the $n$-dimensional Klein bottle as the quotient 
$K_n\simeq (S^1)^n/(z_1,z_2,\hdots ,z_{n-1},t)\sim (\overline{z}_1,\overline{z}_2,\hdots, \overline{z}_{n-1},-t)$ 
and let $\sigma_1:K_n\to K(\ZZ/2,1)$ be the map classifying the double cover. 
\item {\bf (i)} We have the recursive formulas which relate twisted and untwisted KO-groups
\begin{align}
&\widetilde{\op{KO}}^i(K_2)\cong  \widetilde{\op{KO}}^i(\RR P^1)\oplus
  \widetilde{\op{KO}}^i(\RR P^2) \label{base1kn}
\\
&\widetilde{\op{KO}}^i_{\sigma_1}(K_2) \cong \widetilde{\op{KO}}^{i+1}(\RR P^2)\oplus
 \widetilde{\op{KO}}^i (\RR P^2) \label{base2kn}
\\
&\widetilde{\op{KO}}^i(K_{n}) \cong \widetilde{\op{KO}}^i(K_{n-1})\oplus 
\widetilde{\op{KO}}^{i}_{\sigma_1}(K_{n-1}) \label{reckn}
\end{align}
and which, from Table  \ref{korealprfu}, determines the $\op{KO}$-theory of all 
higher-dimensional Klein bottles. 
 
\item {\bf (ii)} We also have the recursive formula for the twisted $\op{KO}$-theory
\(
\widetilde{\op{KO}}_{\sigma_1}^i(K_{n}) \cong \widetilde{\op{KO}}_{\sigma_1}^i(K_{n-1})\oplus 
\widetilde{\op{KO}}^{i}_{\sigma_1}(K_{n-1}) \;.
\label{reckntw}
\)
\end{proposition}
\theproof
{\bf (i)} 
 We begin by verifying the base cases \eqref{base1kn} and \eqref{base2kn}. To compute the $\op{KO}$-theory 
of the 2-dimensional Klein bottle $K_2$, observe  that $K_2$ can be identified as the unit sphere bundle of 
the rank 2 vector bundle 
$$
\xymatrix{
\mathbf{1}\oplus \xi\simeq \RR^2\times S^1/(u,v,t)\sim (u,-v,-t)\ar[r]&\RR P^1
},
$$
where $\xi\to \RR P^1$ is the canonical bundle. This sphere bundle admits a section $s:\RR P^1\to K_2$ 
given by sending  $[x]\mapsto [(1,0,x)]$ with $x\in S^1$. We thus have a homotopy commutative diagram 
\(\label{diagklrps}
\xymatrix{
K_2\simeq S(\mathbf{1}\oplus \xi)\ar[r]^-{j} & D( \mathbf{1}\oplus \xi)\ar[r]\ar[d]_{\pi}^-{\simeq} &
 {\rm Th}(\mathbf{1}\oplus \xi)\simeq S^1\wedge \RR P^2\ar[dl]
\\
&\ar[lu]^-{s} \RR P^1&
}
\)
with the top sequence a cofiber sequence, where  ${\rm Th}$ denotes the Thom space.
This gives rise to a split short exact sequence in $\op{KO}$-theory
$$
\xymatrix{
0 \ar[r] & \widetilde{\op{KO}}^i(\RR P^1)\ar[r]^-{j^*} &
 \widetilde{\op{KO}}^i(K_2) \ar[r] &
  \widetilde{\op{KO}}^{i+1}(S^1\wedge \RR P^2)
  \ar[r] & 0
  },
$$
which gives the identification \eqref{base1kn}.

Now to compute the twisted theory, observe that the twist $\sigma_1$ corresponding to the double cover of $K_2$ 
factors as $\sigma_1:K_2\to \RR P^1\to K(\ZZ/2,1)$, the latter corresponding to $w_1$ of $\xi$. The maps 
in \eqref{diagklrps} are manifestly compatible with twist defined by pullback. Moreover, the usual Thom 
isomorphism in cohomology implies that the twist vanishes on the Thom space. The diagram \eqref{diagklrps} 
thus also gives rise to a split exact sequence in twisted $\widetilde{\op{KO}}$-theory
$$
\xymatrix{
0 \ar[r] &\widetilde{\op{KO}}^*_{\sigma_1}(\RR P^1)\ar[r]^-{j^*}&
 \widetilde{\op{KO}}^*_{\sigma_1}(K_2)\ar[r]&
 \widetilde{\op{KO}}^{*+1}(S^1\wedge \RR P^2)\ar[r]& 0
 }.
 $$
Combining this with the suspension and twisted Thom isomorphism gives the identification \eqref{base2kn}. 

The diagram \eqref{diagklrps} can be extended to higher dimensions as follows. The $n$-dimensional Klein bottle
$K_{n}$ can be identified with the unit sphere bundle of the plane bundle $\RR^2\times (S^1)^{n-1}/\sim \to K_{n-1}$ 
and this bundle again has a section given by $s:[x]\mapsto [(1,0,x)]$. We have a homotopy commutative diagram 
\(\label{diagklrps2}
\xymatrix{
K_n\simeq S(\mathbf{1}\oplus \xi)\ar[r]^-{j} & D( \mathbf{1}\oplus \xi)\ar[r]\ar[d]_{\pi}^-{\simeq} & 
{\rm Th}(\mathbf{1}\oplus \xi)\ar[dl]
\\
&\ar[lu]^-{s}K_{n-1}&
}
\)
with $\xi\to K_{n-1}$ the line bundle classified by $\sigma_1:K_{n-1}\to \RR P^1\to K(\ZZ/2,1)$. This gives a split short exact sequence
$$
\xymatrix{
0 \ar[r]& \widetilde{\op{KO}}^*(K_{n-1})\ar[r]^-{j^*} &
\widetilde{\op{KO}}^*(K_n)\ar[r]&
 \widetilde{\op{KO}}^{*+1}({\rm Th}(\mathbf{1}\oplus \xi))\ar[r]& 0
}.
$$
From the twisted Thom isomorphism, we get the identification \eqref{reckn}. 

\noindent {\bf (ii)} Finally, the recursive formula for the twisted case 
follows again by the diagram \eqref{diagklrps2}, the twisted Thom isomorphism, and the 
compatibility of the twists on $K_{n-1}$ and $K_n$.
\endofproof

Given these preliminary identifications, we return to the question of identifying the remaining coefficients in
 expressions \eqref{dfeqns1}, \eqref{dfeqns2}, and \eqref{dfeqns3} of Remark \ref{rem-pos}.

\begin{lemma}\label{Lemm-cbprime}
In equation \eqref{dfeqns2}, we have $c^{\prime}=0$ and $b^{\prime}=1$.
\end{lemma}
\theproof
Let $\xi\to \RR P^n$ denote the canonical bundle. The Thom space of the Whitney sum $k\xi$ is
 given by the stunted projective spaces
$$
(\RR P^n)^{k\xi}=\RR P^{n+k}/\RR P^{k-1}\;.
$$
For $k=1$, we have $(\RR P^n)^{\xi}=\RR P^{n+1}$. Let $\sigma_1=w_1:\RR P^n\to K(\ZZ/2,1)$ denote the 
generator, corresponding to $w_1(\xi)$ and let $\mathcal{Z}$ denote the local coefficient system 
with fiber $\ZZ$, classified by $w_1$. From the well-known expression (see \cite{GS6} for an explicit calculation) 
$$
H^p(\RR P^{2n}; \Z)=\left\{ \begin{array}{cc}
\ZZ & p={2n}
\\
\ZZ/2 & p\ \text{odd}
\\
0 &\text{otherwise}
\end{array}
\right.
\qquad
\text{and}
\qquad
H^p(\RR P^{2n+1}; \Z)=\left\{ \begin{array}{cc}
\ZZ/2 & p\ \text{odd}
\\
0 &\text{otherwise}
\end{array}
\right.
$$
we identify the relevant terms in the $E_2$-page of the twisted AHSS as 
\( \label{spseqrpn1}
\begin{tikzpicture}[baseline=(current bounding box.center)] 
\matrix (m) [matrix of math nodes,
nodes in empty cells,nodes={minimum width=3ex,
minimum height=3ex,outer sep=0pt},
column sep=1ex,row sep=1ex]{
     &   \phantom{-}0   &  0 & \ZZ/2 & 0 & \ZZ/2 & 0 & \hdots  \\
     &  -1    & \ZZ/2& \ZZ/2 &  \ZZ/2 &  \ZZ/2 &  \ZZ/2 &\\
      & -2   & \ZZ/2& \ZZ/2 &  \ZZ/2 &  \ZZ/2 &  \ZZ/2 & \hdots \\
        & &  &  & \\
     & -4  &  0 & \ZZ/2 & 0 & \ZZ/2 & 0 & \hdots \\
\quad\strut &   \strut \\ };
\draw[thick] (m-1-2.east) -- (m-5-2.east) ;
\end{tikzpicture}
\vspace{-1cm}
\)
where the numbers to the left of the vertical line correspond to the grading on $\pi_q(\op{KO})$ and the 
groups on the right are organized from left to right according to cohomological degree $H^p$ (as usual). 
The entries are zeros for $p>n$ for $n$ odd,  and $ \ZZ$'s are inserted at the entries $E^{-4q,n}_2$ when $n$ is even. 

From the identification of the $\op{KO}$ groups of real projective space in Table \ref{korealprfu} and the twisted 
$\op{KO}$-groups, identified in Remark \ref{twkorpntiso}, we have 
$\op{KO}^0_{w_1(\xi)}(\RR P^3) \cong \widetilde{\op{KO}}^1(\RR P^4)\cong \ZZ/2$. 
Hence, precisely one of the $\ZZ/2$ factors must be killed by either $d_2^{0,-1}$ or $d^{1,-1}$. This is not possible if $b^{\prime}=0$ and $c^{\prime}=1$, since otherwise for the generator $x\in H^1$, we would have $d_2^{0,-1}(x)=d^{1,-1}(x)=x^3\neq 0$ and both factors would be killed. It is also clearly not possible if $b^{\prime}=c^{\prime}=0$. Hence we must have $b^{\prime}=1$. 

To see that $c^{\prime}=0$, we use $\RR P^2$. In this case the twisted Thom isomorphism for the canonical bundle $\xi\to \RR P^2$ gives the identification $\op{KO}^0_{w_1(\xi)}(\RR P^2)\cong \widetilde{\op{KO}}^1(\RR P^3)\cong \ZZ/2\oplus \ZZ/2$. The AHSS looks as follows
\( \label{spseqrpn1p}
\begin{tikzpicture}[baseline=(current bounding box.center)] 
\matrix (m) [matrix of math nodes,
nodes in empty cells,nodes={minimum width=3ex,
minimum height=3ex,outer sep=0pt},
column sep=1ex,row sep=1ex]{
     &   \phantom{-}0   &  0 & \ZZ/2 & \ZZ & 0 & \hdots  \\
     &  -1    & \ZZ/2& \ZZ/2 &  \ZZ/2 &  0 &\\
      & -2   & \ZZ/2& \ZZ/2 &  \ZZ/2 &  0& \hdots \\
        & &  &  & \\
     & -4  &  0 & \ZZ/2 & \ZZ & 0 & \hdots \\
\quad\strut &   \strut \\ };
\draw[thick] (m-1-2.east) -- (m-5-2.east) ;
\end{tikzpicture}
\vspace{-1cm}
\)
Hence both $\ZZ/2$'s must contribute. Therefore, $d_2^{0,-1}$ must vanish and, since both $Sq^2$ and $Sq^1$ vanish for degree reasons, this forces $c^{\prime}=0$.
\endofproof


\begin{lemma}
In equations \eqref{dfeqns1}, \eqref{dfeqns2} and \eqref{dfeqns3}, we have $d^{\prime},d^{\prime\prime},d=1$.
\end{lemma}
\theproof
Consider the twist $w_2:\RR P^5\to K(\ZZ/2,2)$ corresponding to the second Stiefel-Whitney class of the tangent bundle. Note that $\RR P^5$ is oriented, so $w_1=0$. Since $T\RR P^5\oplus {\bf 1}\cong 6\xi$, we have the relation
$$
(\RR P^5)^{T\RR P^5}\wedge S^1\simeq (\RR P^5)^{6\xi}\simeq \RR P^{11}/\RR P^5\;,
$$
 which gives an isomorphism $\widetilde{\op{KO}}^{q-1}((\RR P^5)^{T\RR P^5})\cong \widetilde{\op{KO}}^q(\RR P^{11}/\RR P^5)$. 
 From the long exact sequence for $\op{KO}$-theory of a pair, we have the sequence
\(\label{longseqrp2}
\xymatrix@C=1.7em{
\cdots \ar[r]& \widetilde{\op{KO}}^{q-1}(\RR P^5)\ar[r]& 
\widetilde{\op{KO}}^q(\RR P^{11}/\RR P^5)\ar[r] & 
\widetilde{\op{KO}}^q(\RR P^{11})\ar[r]& \widetilde{\op{KO}}^{q}(\RR P^5)\ar[r] & \cdots}.
\)
In particular, the identifications in Table \ref{korealprfu} and the twisted Thom isomorphism (applied to the tangent bundle) 
give the exact sequence
\(\label{exstw11t5}
\xymatrix{
0\ar[r] &  \op{KO}_{w_2}^{-2}(\RR P^5)\ar[r] & \ZZ/(2^4)\ar[r]^-{\epsilon}&
 \ZZ/(2^2)\ar[r] &  \op{KO}_{w_2}^{-1}(\RR P^5)\ar[r] & 0
}.
\)
We now examine this sequence in conjunction with the AHSS. Since the degree 1 twist $w_1=0$, we need not consider 
local coefficients in the twisted AHSS and in this case the relevant portion of the spectral sequence reduces to 
\(\label{spseqrpn2}
\begin{tikzpicture}[baseline=(current bounding box.center)] 
\matrix (m) [matrix of math nodes,
nodes in empty cells,nodes={minimum width=3ex,
minimum height=3ex,outer sep=0pt},
column sep=1ex,row sep=1ex]{
     &   \phantom{-}0   &   \ZZ & 0 & \ZZ/2 & 0 & \ZZ/2 & \ZZ & \hdots  \\
     &  -1    & \ZZ/2& \ZZ/2 &  \ZZ/2 &  \ZZ/2 &  \ZZ/2 & \ZZ/2 & \\
      & -2   & \ZZ/2& \ZZ/2 &  \ZZ/2 &  \ZZ/2 &  \ZZ/2 & \ZZ/2 & \hdots \\
        & &  &  & \\
     & -4    &   \ZZ & 0 & \ZZ/2 &  0& \ZZ/2 & \ZZ & \hdots \\
\quad\strut &   \strut \\ };
\draw[thick] (m-1-2.east) -- (m-5-2.east) ;
\end{tikzpicture}
\vspace{-1cm}
\)
From \eqref{spseqrpn2}, one sees that both copies of $\ZZ/2$ must contribute to $\widetilde{\op{KO}}_{\sigma_2}^{-2}(\RR P^5)$ and hence this group has 4 elements. Considering the possible images of $\epsilon$ in \eqref{exstw11t5}, the first isomorphism theorem implies that this is only possible if $\epsilon$ is surjective. This, then implies that $\widetilde{\op{KO}}^{-1}_{\sigma_2}(\RR P^5)\cong 0$ and this is only possible if $d_2^{0,-1}$ and $d_3^{1,-2}$ are isomorphisms. Since $Sq^2$ vanishes here for degree reasons and $w_1=0$, we must have $d^{\prime\prime}=1$ and $d^{\prime}=1$.

To show that $d=1$, we note that the sequence \eqref{longseqrp2}, the identifications in Table \eqref{korealprfu}, 
and the twisted Thom isomorphism give an exact sequence
\(\label{sesporrp5}
\widetilde{\op{KO}}_{w_2}^{3}(\RR P^5)\longrightarrow \ZZ/2\oplus \ZZ/2\overset{\phi}{\longrightarrow} \ZZ 
\)
and $\phi$ is necessarily trivial. This gives a surjection $\widetilde{\op{KO}}_{w_2}^{3}(\RR P^5)\onto \ZZ/2\oplus \ZZ/2$. 
From the AHSS it then follows that on the generator $x^2\in H^2$, we must have $d^{2,0}_2(x^2)=x^4+dx^4=0$. Hence $d=1$. 
\endofproof

The final identification requires the Klein bottles and dispenses of the $w_1^2$ term. 

\begin{lemma}
In equation \eqref{dfeqns1}, we have $b=1$ and $c=0$. 
\end{lemma}
\theproof
First, applying the AHSS to $\RR P^3$ as in the proof of Lemma \ref{Lemm-cbprime}, we see that the differential 
$d^{1,0}_2:\ZZ/2\to \ZZ/2$ must be an isomorphism. Hence, either $b=0$ and $c=1$ or $b=1$ and $c=0$. 
We claim it is the latter. Let  $\sigma_1:K_3\to K(\ZZ/2,1)$ be the twist corresponding to the double cover of the 3-dimensional 
Klein bottle $K_3$. This twist satisfies $\sigma_1\neq 0$, but $\sigma_1^2=0$. Moreover Let $\mathcal{Z}$ be the local 
coefficient system with fiber $\ZZ$ corresponding to $\sigma_1$. A straightforward calculation with the Mayer-Vietoris 
sequence in twisted cohomology 
\footnote{One can use the covering of $K_3$ by $(S^1)^2\times (0,1)$ and $K_3-((S^1)^2\times\{\frac{1}{2}\})$.} 
shows that $H^1(K_3; \mathcal{Z})\cong \ZZ/2$. On the other hand, our calculation of the twisted $\op{KO}$-theory 
of the Klein bottles gives
$$
\widetilde{\op{KO}}^{1}_{\sigma_1}(K_3)\cong 
\widetilde{\op{KO}}^{1}_{\sigma_1}(K_2)\oplus \widetilde{\op{KO}}^{1}_{\sigma_1}(K_2)\cong (\ZZ/2)^4\;.
$$
Combining this with the identification 
$
H^i(K_3; \Z/2)=(\ZZ/2)^{3\choose i}
$,
obtained in \cite{Dav} and our identification of the differential $d_2^{p,8t-1}$, it follows that the differential 
$d_2^{1,0}:H^1(K_3; \mathcal{Z})\to \ZZ/2\cong H^3(K_3;\ZZ/2)\cong \ZZ/2$ must be an isomorphism.
Hence $b=1$ and $c=0$. 
\endofproof

Assembling the above lemmas together, we have proved the following.

\begin{proposition}[AHSS for twisted KO-theory.]\label{prop-fullAhss}
 Let $X$ be a CW complex. Denote the full twisting  class for KO-theory by 
$\sigma=(\sigma_1, \sigma_2) \in  H^1(X; \ZZ/2) \times H^2(X; \ZZ/2)$.
Then  there exists a spectral sequence for twisted KO-theory with
$E_2$-term the twisted {\v C}ech cohomology 
$$
E_2^{p, q}= \check{H}^p(X, \pi^{q}(KO)_{\sigma_1})\;,
$$
with differentials $d_2^{p, -8t}: \check{H}^p(X; \ZZ) \longrightarrow \check{H}^{p+2}(X; \ZZ/2)$, 
$d_2^{p, -8t-1}: \check{H}^p(X; \ZZ/2) \longrightarrow \check{H}^{p+2}(X; \ZZ/2)$,
and 
$d_3^{p, -8t-2}: \check{H}^p(X; \ZZ/2) \longrightarrow \check{H}^{p+3}(X; \ZZ)$
given, respectively, by
\bea
d_2^{p, -8t}&=& Sq^2 \circ r + \sigma_1 \cup Sq^1 \circ r+\sigma_2\cup r \;,
\\
d_2^{p, -8t-1}&=& Sq^2   +  \sigma_1 \cup Sq^1 +\sigma_2\cup \;,
\\
d_3^{p, -8t-2}&=& \beta \circ Sq^2  + \beta\circ \sigma_2\cup\;.
\eea
\end{proposition}

Note the modification of the $E_2$-page upon including the degree one twist. 

\subsection{AHSS for twisted differential KO-theory}
\label{Sec-AHSS-hatKO}
\begin{figure}[h!]\label{theE2sps}
\hspace{-4mm}
\begin{minipage}{.5\textwidth}
\centering
\footnotesize
\begin{tikzpicture}
\matrix (m) [matrix of math nodes,
nodes in empty cells,nodes={minimum width=3ex,
minimum height=3ex,outer sep=0pt},
column sep=1ex,row sep=1ex]{
         \phantom{-}8    &  H^0(M;\ZZ) &  H^1(M;\ZZ)  & H^2(M;\ZZ)   \\
          \phantom{-} 7  & H^0(M;\ZZ/2)& H^1(M;\ZZ/2) &  H^2(M;\ZZ/2) \\
           \phantom{-}6  & H^0(M;\ZZ/2)& H^1(M;\ZZ/2) &  H^2(M;\ZZ/2) \\
             & & & \\
        \phantom{-} 4    &  H^0(M;\ZZ) &  H^1(M;\ZZ)  & H^2(M;\ZZ)   \\
             & & &  \\
             & & & \\
             & & &  \\
          \phantom{-}0   &  \Omega^0_{\rm cl,\ZZ}(M;\pi_*(KO))  & & \\
      -1    & H^0(M;\ZZ/2)& H^1(M;\ZZ/2) &  H^2(M;\ZZ/2)\\
       -2   & H^0(M;\ZZ/2)& H^1(M;\ZZ/2) &  H^2(M;\ZZ/2) \\
      -3.   & H^0(M;U(1)) &  H^1(M;U(1))  & H^2(M;U(1))  \\
             & & &\\
             & & &  \\
             & & & \\
        -7  & H^0(M;U(1)) & H^1(M;U(1))  &  H^2(M;U(1))  \\
\quad\strut &   \strut \\ };
\draw[thick] (m-1-1.east) -- (m-16-1.east) ;
\end{tikzpicture}
\end{minipage}%
~~~~~~~~~~~~~
\begin{minipage}{.5\textwidth}
\centering
\footnotesize
\begin{tikzpicture}
\matrix (m) [matrix of math nodes,
nodes in empty cells,nodes={minimum width=3ex,
minimum height=3ex,outer sep=0pt},
column sep=1ex,row sep=1ex]{
        \phantom{-} 8   &  H^0(M;\ZZ) &  H^1(M;\ZZ)  & H^2(M;\ZZ)   \\
           \phantom{-}7 & H^0(M;\ZZ/2)& H^1(M;\ZZ/2) &  H^2(M;\ZZ/2) \\
          \phantom{-} 6 & H^0(M;\ZZ/2)& H^1(M;\ZZ/2) &  H^2(M;\ZZ/2) \\
            & & &\\
         \phantom{-}4   &  H^0(M;\ZZ) &  H^1(M;\ZZ)  & H^2(M;\ZZ)   \\
            & & & \\
            &  & &\\
            &  & & \\
 \phantom{-}0  & H^0(M;\ZZ) & H^1(M;\ZZ) & H^2(M;\ZZ)  \\
       -1  & H^0(M;\ZZ/2)& H^1(M;\ZZ/2) &  H^2(M;\ZZ/2)\\
       -2  & H^0(M;\ZZ/2)& H^1(M;\ZZ/2) &  H^2(M;\ZZ/2) \\
            & & &  \\
        -4 &H^0(M;\ZZ) &  H^1(M;\ZZ)  & H^2(M;\ZZ) \\
            & & & \\
            & & & \\
            & & &  \\
        -8 &H^0(M;\ZZ) &  H^1(M;\ZZ)  & H^2(M;\ZZ) \\        
\quad\strut &   \strut \\ };
\draw[thick] (m-1-1.east) -- (m-17-1.east) ;
\end{tikzpicture}
\end{minipage}
\vspace{-10mm}
\caption{The $E_2$-page for $\widehat{\rm KO}$ on the left vs. the $E_2$-page for $\op{KO}$ on the right.}
\end{figure}
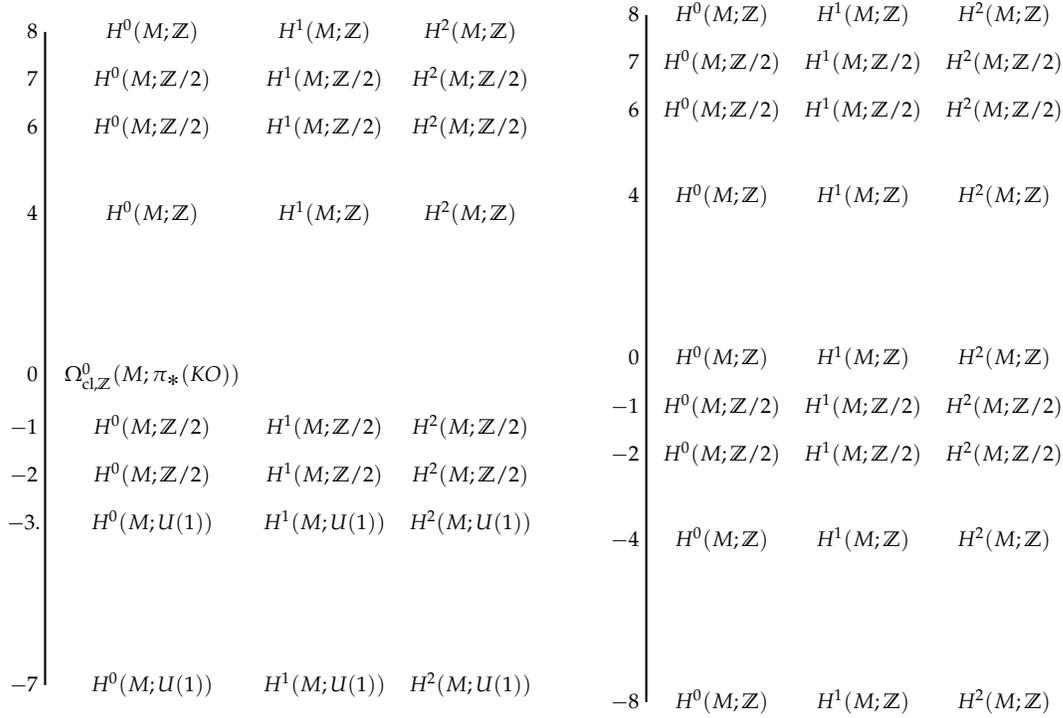

With the heavy lifting done, we now move to the differential theory. 
For twists of the form $\sigma_2 : M\to \BB^2\ZZ/2$, and in contrast to the case of twists by $\sigma_1$
in Proposition \ref{prop-fullAhss}, the spectral sequence in Proposition \ref{spseqdfko} 
has $E_2$-page which looks identical to the untwisted case. \footnote{Note that in general this is not the case, since for 
non-torsion twists the curvature forms in bidegree $(0,0)$ may be twisted closed.} 
In \cite{GS-KO}, we identified the untwisted $E_2$-page as in figure \eqref{theE2sps}.

If we include the degree 1 differential twists $\hat{\sigma}_1:M\to \BB(\ZZ/2)_{\nabla}$, then the topological AHSS is 
modified to include cohomology with local coefficients. In contrast, the differential AHSS has both local coefficients 
as well as the modification 
$$\xymatrix{
E^{0,0}_2=\Omega_{\rm cl,\ZZ}^0(M;\pi_*(\op{KO}))
\ar@{~>}[r] &
 \; E^{0,0}_2=\Omega^0(M;\mathcal{L})_{{\rm fl}\text{-}\nabla}
 },
$$
where $\mathcal{L}\to M$ is the 4-periodic graded bundle with fiber $\RR[\alpha,\alpha^{-1}]$, $\vert \alpha\vert=4$, 
with classifying map given by the topological twist $\sigma_1$ and flat connection $\nabla$ specified by the 1-form 
$\mathcal{R}(\hat{\sigma}_1)$ (see Proposition \ref{prop-prop}). 

\medskip
We now identify a notable portion of the differentials in the twisted AHSS.  As in \cite{GS-KO}, we split this task into 
two parts -- the identification of the flat differentials and the identification of the geometric differentials. 
We begin with the former.

\begin{proposition}[Flat differentials]\label{fldfko}
Let $\hat\sigma=(\hat\sigma_1, \sigma_2) \in  H^1(X; \ZZ/2) \times H^2(X; \ZZ/2)$ be a twist of $\op{KO}$-theory 
and regard this as a twist for the  differential refinement via Proposition \ref{dftwtptw}. Let $\mathcal{L}^{\delta}$ 
and $\mathcal{Z}$ be the local systems with fiber $\RR$ and $\ZZ\subset \RR$, respectively, both classified by $\sigma_1$. 
Then the differentials 
\begin{align*}
& d_3^{p, -8t-7}: \check{H}^p(X; \mathcal{L}^{\delta}/\mathcal{Z}) \longrightarrow \check{H}^{p+3}(X; \ZZ/2), \nonumber
\\
& d_2^{p, -8t-1}: \check{H}^p(X; \ZZ/2) \longrightarrow \check{H}^{p+2}(X; \ZZ/2), \nonumber
\\
& d_3^{p, -8t-2}: \check{H}^p(X; \ZZ/2) \longrightarrow \check{H}^{p+2}(X; \mathcal{L}^{\delta}/\mathcal{Z})
\end{align*}
are given, respectively, by
\bea
d_2^{p, -8t-7}&=& Sq^2 \circ r \circ \beta + \sigma_1 \cup Sq^1 \circ r \circ \beta +\sigma_2\cup r\circ \beta\;,
\\
d_2^{p, -8t-1}&=& Sq^2   + \sigma_1  \cup Sq^1 + \sigma_2 \cup \;,
\\
d_3^{p, -8t-2}&=& j \circ Sq^2  + j\circ (\sigma_2 \cup)\;,
\eea
where $\beta:H^p(M;\mathcal{L}^{\delta}/\mathcal{Z})\to H^{p+1}(M;\mathcal{Z})$ is the Bockstein 
associated to the exact sequence of locally constant sheaves
$
0\longrightarrow \mathcal{Z}\longrightarrow \mathcal{L}^{\delta}\longrightarrow \mathcal{L}^{\delta}/\mathcal{Z}
\longrightarrow 0\;.
$
\end{proposition}
\theproof
From the proof of \cite[Proposition 17]{GS3}, \footnote{Note that $\op{KO}_{\sigma_1}$ does not fulfill the requirements 
of the referenced 
proposition in \cite{GS3}. However, the proposition still holds at the entries where the connecting homomorphism 
of spectral sequences is trivial and the homological argument presented there works equally well for twisted cohomology 
and local coefficients.} 
the Bockstein for cohomology with local coefficients $\beta:H^p(M;\mathcal{L}^{\delta}/\mathcal{Z})\to H^{p+1}(M;\mathcal{Z})$ 
commutes with the differentials in the spectral sequences 
which appear in the fourth quadrant $p>0,q<0$, corresponding to $\widehat{\op{KO}}_{\hat{\sigma}}$ and 
$\op{KO}_{\sigma}$, respectively. 
Consideration of the morphisms of AHSS's induced by the long fiber sequence of parametrized spectra on $M$, associated to 
$$
\op{KO}_{\sigma}\longrightarrow \op{KO}_{\sigma}\wedge \mathbb{S}(\RR)
\longrightarrow \op{KO}_{\sigma}\wedge \mathbb{S}(\RR/\ZZ)\;,
$$
 shows that the connecting homomorphism 
 $\tilde{\beta}:\op{KO}_{\sigma}\wedge \mathbb{S}(\RR/\ZZ)\to \Sigma \op{KO}_{\sigma}$ 
 induces an isomorphism at the entries with $\ZZ/2$-coefficients. This proves the claim for the differentials $d_2^{p,-8t-1}$. 
 For the remaining differentials, observe that from the above considerations, we have commutative diagrams
\(
\xymatrix{
H^p(X;\ZZ/2)\ar@{=}[r]\ar[d]^-{d_2} & H^p(X;\ZZ/2)\ar[d]^-{d_3}
\\
H^{p+2}(X;\mathcal{L}^{\delta}/\mathcal{Z})\ar[r]^-{\beta} & H^{p+3}(X;\mathcal{Z})
}
\qquad
\text{and}
\qquad 
\xymatrix{
H^{p}(X;\mathcal{L}^{\delta}/\mathcal{Z})\ar[r]^-{\beta}\ar[d]^-{d_3} & H^{p+1}(X;\mathcal{Z})\ar[d]^-{d_2}
\\
H^{p+3}(X;\ZZ/2)\ar@{=}[r] & H^{p+3}(X;\ZZ/2)\;.
}
\)
The right diagram gives the identification of $d_3^{p,-8t-7}$. As for the remaining differential, note that the from the 
relation $\mathcal{I}j=\beta$, with $\beta:H^p(X;\mathcal{L}^{\delta}/\mathcal{Z})\to H^{p+1}(X;\mathcal{Z})$ 
the Bockstein, we have that $d_2^{p,-8t-2}$ is as claimed, modulo possibly an element in the kernel of $\beta$. 
But naturality of $d_2$ and the identification 
$$
H^2(K(\ZZ/2,p);\mathcal{L}^{\delta})\cong H^2(K(\ZZ/2,p);\RR)\cong 0\;,
$$
forces this element to vanish.
\endofproof

Note that flat KO-theory $KO(-; \R/\Z)$  has been used in \cite[Appendix B]{FMS} to establish Pontrjagin 
self-duality of KO-theory.

\medskip
We now turn our attention to the geometric differentials in the AHSS. In \cite{GS-RR}, starting with 
a form $\omega = \omega_0 + \omega_2 + \omega_4 + \cdots$,  we identified the following 
differentials for complex differential K-theory
\begin{align}
&d^{0,0}_2(\omega)=[\omega_2]\!\!\mod \ZZ\;,
\\
&d^{0,0}_4(\omega)=[\omega_4]\!\!\mod \ZZ+j_2Sq^2(x_2)\;, \label{d4dfkth}
\\
&d^{0,0}_6(\omega)=[\omega_6]\!\!\mod \ZZ+\overline{Sq}^2(\omega_4)+\mathcal{P}_3^1(x_2)\;,
\end{align}
where $x_2$ is an integral lift of $[\omega_2]$, $\overline{Sq}^2(\omega_4)$ is some $\ZZ/2$-valued 
cohomology operations, well-defined modulo the image of $d_3$, and $\mathcal{P}_3^1$ denotes the 
first Steenrod power at the prime 3. Here the identification of the differentials $d_4$ and $d_6$ were made 
more complicated due to the presence of nonvanishing $d_2$ and $d_4$, pairwise respectively. 
For $\op{KO}$, the situation is made easier for $d_4$, since there are no degree 2 forms that arise, 
while the differential $d_8$ appears to be quite complicated. We have the following identification.

\begin{proposition}[Geometric differential]\label{prop-geomdiff}
Consider the fourth differential $d_4:E^{0,0}_4\to E^{5,-4}_4$.
\item {\bf (i)} In the AHSS for untwisted 
$\widehat{\op{KO}}$, $d_4$ is given by 
\(\label{4thdfkodf}
d_4(\omega)=\tfrac{1}{2} [\omega_4]\!\!\mod \ZZ\;.
\)
\item {\bf (ii)} For twisted $\widehat{\op{KO}}$, with twist 
$\hat\sigma=(\hat\sigma_1,\sigma_2): M\longrightarrow \BB(\ZZ/2)_{\nabla}\times \BB^2 \ZZ/2$, 
the  differential is given by 
\(\label{4thdfkodf}
d_4(\omega)=\tfrac{1}{2} [\omega_4]\!\!\mod \mathcal{Z}\,,
\)
where $[\omega_4]$ is the class of a flat 4-form, with values in the graded line bundle $\mathcal{L}$, 
classified by the rationalization of $\sigma_1$. The local system $\mathcal{Z}\into \mathcal{L}$ is 
the $\mathbb{Z}$-subbundle of the degree 4 component (with respect to the aforementioned
grading) of $\mathcal{L}$, classified by $\sigma_1$.  
\end{proposition}
\theproof
As in \cite{GS-RR}, we construct a universal example as follows. The differential 
$$
d_4:E^{0,0}_4=\bigoplus_n\Omega^{4n}_{\rm cl}(M)\longrightarrow H^4(M;U(1))/{\rm im}(d_2)
$$ 
factors through the cohomology group $H^4(M;\RR)$. We can approximate the Eilenberg MacLane space $K(\RR,4)$ by 
a smooth manifold via a sequence of surgeries and the differential $d_2=Sq^2$ necessarily vanishes on $K(\RR,4)$. 
There is no cokernel, and it follows that $d_4$ must be a cohomology operation of the form $\phi:K(\RR,4)\to K(U(1),4)$. 
Such operations are parametrized by $\lambda\in \RR$, corresponding to the map on cohomology 
$x\mapsto \lambda x\!\!\mod \ZZ$.  To see that $\lambda=\frac{1}{2}$, we consider the 4-sphere $S^4$. 
The Pontrjagin character maps $\alpha\in \op{KO}^{-4}(\ast)\mapsto 2\in \ZZ$ and hence a form 
$\omega_4\in  \Omega_{\rm cl}(S^4)$ is identified with an element of the image of ${\rm Ph}$ if and only 
if it has even integral periods -- hence $\frac{1}{2}[\omega_4]$ must be integral. 

The twisted case is similar. In this case, $d_4$ factors through $H^4(M;\cL^\delta)$ and the twisted cohomology is 
represented by sections of the pullback of the bundle $K(\RR,4)\sslash\RR^{\times}\to K(\RR^{\times},1)$, 
where the units $\RR^{\times}$ act by multiplication on $\RR$ in the usual way. The maps $\phi^\prime:K(\RR,4)\to K(\RR,4)$ 
defined by multiplication by $\lambda$ are equivariant with respect to the $\RR^{\times}$-action. Moreover, all the maps
$$
\phi:K(\RR,4)\!\sslash\!\RR^{\times}\longrightarrow
 (K(\RR,4)\!\sslash\!\RR^{\times})\sslash(K(\ZZ,4)\!\sslash\! \ZZ/2)
$$
factor, up to homotopy, through a map 
$$
\phi^{\prime}:K(\RR,4)\!\sslash\!\RR^{\times}\longrightarrow
 K(\RR,4)\!\sslash\!\RR^{\times}\;.
$$
The homotopy classes of $\phi^{\prime}$, and hence $\phi$ are thus in bijective correspondence with 
$\lambda\in \RR$, as in the untwisted case. Since we have already shown that $\lambda=\frac{1}{2}$ 
in the untwisted case, the claim follows.
\endofproof

\subsection{Applications}
\label{Sec-apps} 

Having provided our constructions in previous sections, we now illustrate with applications 
originating from physics, but which will interestingly lead to results in geometry and topology,
 including obtaining Rokhlin's theorem as a corollary. 

\subsubsection*{\bf (I) Low-dimensional manifolds} 
For low-dimensional complexes, there is an isomorphism  of the form 
$$
KO(X) \cong \Z \times H^1(X; \Z/2) \times H^2(X; \Z/2)\;.
$$
This is classical when $X$ is a surface, by  Atiyah \cite{At-s} (using additive notation)
and has been used in \cite{MMS} and more recently by Hitchin in \cite{Hi}. 
For $\dim{X} \leq 3$, this is done recently (using product notation)
by Karoubi, Schlichting, and Weibel in
 \cite[proof of theorem 4.6]{KSW} in the context of real varieties. 
  In \cite[Remark 5.9.1]{KW2}, the AHSS is used to show the following: 
  for $\dim X \leq 3$, the kernel of  
  $(\rm rank, \det): KO(X) \to H^0(X; \ZZ)\oplus H^1(X; \ZZ/2)$ 
  is mapped isomorphically onto $H^2(X; \Z/2)$ via 
    the second Stiefel Whitney class $w_2$.

 \medskip
 We generalize the above to the twisted differential setting. 
 \begin{proposition}[Differential KO-theory of low-dimensional manifolds]
 If $M$ is a manifold of dimension $\leq 3$, then we have a bijection of sets
 $$
 \widehat{\op{KO}}^0(M)\cong \ZZ\times H^1(M;\ZZ/2)\times H^2(M;\ZZ/2)\;,
 $$
 and the canonical map $\mathcal{I}:\widehat{\op{KO}}(M)\to \op{KO}(M)$ is an isomorphism of groups.
 \end{proposition}
 \theproof
The first claim follows immediately from consideration of the differential AHSS, combined with the low 
dimension of $M$. To see that $\mathcal{I}$ induces an isomorphism of groups, one simply observes 
that $\mathcal{I}$ induces an isomorphism of corresponding AHSS. Hence, $\widehat{\op{KO}}(M)$ 
fits into the same extension as $\op{KO}(M)$.
 \endofproof
 
 Note that the differential refinement is not reflected in having differential cohomology
 on the right hand side of the isomorphism, due to the dimensions considered versus those
 related to KO. We will see similar effects below, but we emphasize that in higher dimensions
 one would encounter twisted differential integral cohomology (i.e., twisted Deligne cohomology)
 developed in \cite{GS4}\cite{GS6}, but  we will not develop this here.

\medskip
\subsubsection*{\bf (II) Lifting forms to $\widehat{\rm KO}_{\widehat{\rm tw}}$ via field quantization in string theory.}
In the presence of an H-flux or H-field, charges of Type I D-branes  \cite{Wi} and corresponding 
Ramond-Ramond (RR-) fields  \cite{MW} are proposed to take values in twisted KO-theory.  
We have put on firm mathematical ground the statement in the untwisted case \cite{GS-KO} 
and the twisted complex case \cite{GS-RR}; here we do something analogous in the twisted real case. 
We explore the RR-field quantization (borrowed from the type II case \cite{GS-RR}) 
by utilizing the differential AHSS. Although it is possible 
to find the complete quantization conditions for the RR-fields on a 10-dimensional spacetime $X^{10}$, 
this would be somewhat involved and requires identifying higher differentials in the AHSS. The pleasant 
outcome is that what we do already leads to interesting topological and geometric consequences. 
We start with the untwisted case. 

\begin{proposition}[Lifting $4k$-forms to differential KO-theory]\label{Prop-4k}
We have the following necessary and sufficient conditions on a differential form $G=G_0+G_4$ in order 
for it to lift to $\widehat{{\rm KO}}$-theory:
\item {\bf (i)} $G_0$ is a constant integer.
\item {\bf (ii)} $G_4$ satisfies the equation 
$$
[G_4] \; ({\rm mod}\; \ZZ)=j_2Sq^2(x_2)=j_2x_2^2\;,
$$
for some $x_2\in H^2(X;\ZZ/2)$, where $j_2:H^2(M;\ZZ/2)\to H^2(M;U(1))$ is the inclusion
 via the 2-roots of unity.
\end{proposition}
\theproof
Applying the differential AHSS, we see that the only $d_2=j_2Sq^2:H^2(X;\ZZ/2)\to H^4(X;U(1))$ and 
$d_3:H^1(X;\ZZ/2)\to H^4(X;U(1))$ can affect the values of $G$. The differential $d_3$ can be shown 
to vanish for degree reasons, and hence the only obstruction is $j_2Sq^2$; the corresponding vanishing 
condition on $d_4(G_4)=\frac{1}{2}[G_4] \!\!\mod \ZZ$ is thus the one claimed.
\endofproof

\begin{remark}[Rohklin's Theorem]
An immediate corollary of Prop. \ref{Prop-4k} is Rokhlin's Theorem \cite{Rok}. Indeed, consider a 4-dimensional 
Spin manifold $X$. We take $G_4$ to be the degree 4 component of the $\widehat{A}$-genus (the Chern character of 
the virtual spinor bundle) and note that the Wu formula and the fact that $X$ is Spin imply $Sq^2=0$ on 
degree 2 classes. By Hirzebruch's signature formula, we see that the divisibility of $G_4$ by 2 implies divisibility 
of the signature by 16.
\end{remark}

We now turn to the twisted case. Here the torsion twist has an affect on the possible values of $G$, even in the case of 
vanishing degree 1 twist (in which case twisted closed forms reduce to untwisted ones). 

\begin{proposition}[Lifting $4k$-forms to twisted differential KO-theory]\label{lift-formKO}
Let $\sigma=\sigma_2$ be a twist of degree 2 on $X^{10}$ (i.e., the $B$-field). We have the following necessary and 
sufficient conditions on a differential form $G=G_0+G_4$ in order that it lift to $\widehat{\op{KO}}_{\sigma_2}$:
\item {\bf (i)} $G_0$ is integral.
\item {\bf (ii)} We have 
$$
\tfrac{1}{2}[G_4]\; ({\rm mod}\; \ZZ)= j_2x^2_2+j_2\sigma_2\cup x_2\;,
$$
for some $x_2\in H^2(X;\ZZ/2)$. 
 \end{proposition}
\theproof
The proof is completely analogous to Proposition \ref{Prop-4k}. We simply utilize the identification of the differential 
$d_2$ in the twisted differential case from Proposition \ref{fldfko}.
\endofproof

This result in the twisted case also has an interesting consequence for orientable 4-manifolds. 
Taking $\sigma_2=w_2$, the Wu formula implies 
$$
x_2^2+w_2x_2=2x_2^2=0\;.
$$
Hence, regardless of whether the manifold is spinnable, it is necessary that $G_4$ has even integral periods in order 
for it to lift to $\op{KO}$. It is not true that on a general orientable 4-manifold that the signature is divisible by 16. 
This does not contradict the previous result, however, since for non Spin manifolds, we cannot take $G_4$ to be 
the degree 4 component of the $\widehat{A}$-genus (i.e., we do not have a spinor bundle available). 
We have the following:

\begin{corollary}[Periods for 4-manifolds]
For the twist $\sigma_2=w_2(TM)$, where $M$ is an orientable 4-manifold, it is necessary that a differential form 
$G_4$ on $M$ must have even integral periods in order for it to lift to twisted differential $\op{KO}$-theory.
\end{corollary} 

Note that this killing of $Sq^2$ is a typical effect in twisted cohomology and is really a manifestation 
of the twisted Thom isomorphism (see also Remark \ref{twkorpntiso} and Remark \ref{rem-untw}).

\subsubsection*{\bf (III) Twisted differential KO-theory and type I string theory} 

\medskip
Type I string theory with non-trivial B-field is often referred to as type I theory
"without vector structure" with characteristic class $w\in H^2(X, \Z/2)$ for ${\rm PO}(\mathcal{H})\simeq K(\Z/2, 2)$
bundles \cite{Wi-vect}. 
 Twistings distinguishing various forms of  type I string theory as well as corresponding dualities 
are studied in \cite{DMDR}\cite{DMDR2}.
 The corresponding RR-fields may be modeled as principal ${\rm Spin}(32)$-bundles 
 "without vector structure" \cite{LMST} \cite{Wi-vect}.
 The B-field in Type I is a local closed 2-form 
on $M$ which defines a class $B$ in $H^2(M; \Z/2)$ \cite{SeSe}.
A geometric interpretation of such 2-torsion B-fields as the holonomy of connections 
for real bundle gerbes are given in \cite{MMS}, where isomorphism classes of real bundles 
gerbes are shown to be  in bijective correspondence with  Dixmier-Douady
class in $H^2(M; \Z/2)$, known in the physics literature as the
t'Hooft class. 
The anomaly cancellation condition is
$
\widetilde{w}_2(V) - B = 0 \in H^2(M; \Z/2)
$,
where $V\to M$ is an ${\rm SO}(32)/(\ZZ/2)$-bundle "without vector structure" and $\widetilde{w}_2(V)$ 
is the obstruction class for true ${\rm SO}(32)/(\ZZ/2)$-structure. 


\medskip
In order to be self-contained and to show explicitly 
what happens in the differential refinement, we derive this anomaly cancellation from the differential 
refinement of the Freed-Witten anomaly, discussed in \cite{GS-RR}. 
Type I theories without vector structure and corresponding anomalies are discussed in the 
context of twisted KO-theory in \cite{Fr}. What we do here is streamline the description 
to highlight direct connection to $\widehat{\rm KO}_{\widehat{\rm tw}}$, with consequences below in {\bf (IV)}.

\medskip
Let $Q\into M$ be a $D$-brane. The class $w_2(Q)$ gives rise to a differential cohomology class 
$\widehat{W}_3(Q)\in \widehat{H}^3(Q;\ZZ)$ via the map 
$$
j\circ j_2:H^2(Q;\ZZ/2)\longrightarrow H^2(Q;U(1))\longrightarrow \widehat{H}^3(Q;\ZZ)\;,
$$
where $j:H^2(Q;U(1))\to \widehat{H}^3(Q;\ZZ)$ is the inclusion of flat 
classes. This indeed defines a differential refinement of $W_3$, but it is not unique. However, the $U(1)$-gauge field 
$\mathcal{A}$ determines a differential refinement of $W_3$ (see \cite{GS-RR}), where we put forth a differential 
refinement of the Freed-Witten anomaly. Explicitly, this is given as
\(\label{fwatotype1}
\widehat{W}_3(Q;\mathcal{F})+j(B)=j\circ j_2w_2(Q)-j{\rm exp}\big(\tfrac{1}{2\pi i}\mathcal{F}\big)+\hat{h}=0\;,
\)
where $\hat{h}$ is the degree three gerbe and 
$\mathcal{F}$ represents the "curvature" of the $U(1)$-gauge field. Applying this to a spacetime
 filling 9-brane in the orientifold  setting, 
 \footnote{For our purposes it  suffices to view an orientifold as essentially a quotient space of an involution.}  
 we note that since $M$ is ${\rm Spin}$ we have the condition $w_2(Q)=0$,
 and so \eqref{fwatotype1} reduces to 
%
\(\label{bfildcntpe1}
j{\rm exp}\big(\tfrac{1}{2\pi i}\mathcal{F}\big)=\hat{h}\;.
\)
Now the class ${\rm exp}\big(\tfrac{1}{2\pi i}\mathcal{F}\big)\in H^2(Q;U(1))$ is 2-torsion, since $\frac{1}{2\pi i}\mathcal{F}$ 
has half-integral periods (i.e., $\mathcal{F}$ is not the curvature of a true $U(1)$-gauge bundle in general) and hence in 
$\hat{h}=j\circ j_2(B)$ for some class $B\in H^2(M;\ZZ/2)$.
%

\medskip
Now the spacetime filling branes support ${\rm SO}(32)/(\ZZ/2)$-bundles $V\to M$ "without vector structure" and 
the obstruction to lifting to a legitimate bundle is the mod 2-class $\widetilde{w}_2(V)$. The hypothesis that the 
background field $B$ twists  the spacetime geometry (including the gauge bundles on the brane) \cite{Wi}
then yields precisely  the following cancellation condition.

\begin{remark}[Type I anomaly cancellation in twisted differential KO-theory]
On an orientifold spacetime, the $\hat{h}$-field necessarily satisfies the condition \eqref{bfildcntpe1} and hence is 
of the form $\hat{h}=j\circ j_2(B)$ for some $B\in H^2(M;\ZZ/2)$. This class is unique modulo the image of 
the mod 2 reduction $\rho_2:H^2(M;\ZZ)\to H^2(M;\ZZ/2)$ on torsion classes in $H^2(M;\ZZ)$, and 
we have the natural cancellation condition for type I string theory
\(
\label{w2tilde}
\widetilde{w}_2(V)- B=0\in H^2(M;\ZZ/2)\;,
\) 
where $V\to M$ is an ${\rm SO}(32)/(\ZZ/2)$-bundle "without vector structure".
\end{remark}

\begin{remark}[Interpretation and consequences] The above discussion leads to the following:
\item {\bf (i)}  
This is one of the instances 
where differential refinement of torsion does not add geometric information, in the sense that expression 
\eqref{w2tilde} is not refined (no hats); cf. Remark \ref{Rem-twdeg2} and Prop. \ref {dftwtptw}; see also below. 
\item {\bf (ii)}  However, we highlight that the geometry pins down a preferred choice of differential refinement! 
\end{remark} 

\subsubsection*{\bf (IV) Orientations and twisted differential Spin structures.} 
Spin structures and their variants play an important role in string theory and M-theory
(see \cite{pf} for a survey). 
Like orientation with respect to KO-theory is given by Spin structures \cite{ABS}, 
orientation with respect to twisted KO-theory is given by twisted Spin structures \cite{Wang}, in the
sense of \cite{MuS}\cite{To}\cite{Dou}.
We now interpret twisted differential Spin structures.
Unlike the case of ${\rm Spin}^c$ structures, we have the following consequence which is 
reminiscent of the fact that the rational cohomology of $B{\rm O}$ vanishes in degrees 1 and 2,
 and hence Chern-Weil forms do not contribute to the differential characteristic classes in the Whitehead
  tower (see \cite{Cech}\cite{9}) -- only the topological torsion information is retained. 
%

\begin{proposition}[Connections vs. torsion] 
Let $\mathcal{I}:\BB{\rm O}_{\nabla}\to \BB {\rm O}$ be the canonical map which forgets the connection. 
\item {\bf (i)} Then we have isomorphisms
\begin{align*} 
\mathcal{I}^*&:\ZZ/2(w_1)\cong H^1(\BB{\rm O};\ZZ/2)\cong H^1(\BB {\rm O}_{\nabla};\ZZ/2),
\\ \mathcal{I}^*&:\ZZ/2(w_2)\cong H^1(\BB{\rm O};\ZZ/2)\cong H^2(\BB {\rm O}_{\nabla};\ZZ/2)\;,
\end{align*}
where $w_1$ and $w_2$ are the first and second Stiefel-Whitney classes, generating the respective 
copies of $\ZZ/2$. 
\item {\bf (ii)} Moreover, we have a pullback diagram of classifying stacks
$$
\xymatrix@C=3.8em@R=1.5em{
\BB {\rm Spin}_\nabla \ar[r] \ar[d] & \ast \ar[d] 
\\
\BB {\rm O}_\nabla \ar[r]^-{(w_1,\;w_2)} & \BB \ZZ/2\times \BB^2\Z/2\;.
}
$$
\end{proposition}
\theproof
The first claim follows from the cohesive adjunction $(\delta^{\dagger}\dashv \Gamma \dashv \delta\dashv \Pi)$ 
on smooth stacks as developed by Schreiber in \cite{Urs} (see also \cite[Proposition 8]{GS3} for how this 
applies in this context). More precisely, we have the following natural isomorphisms
\begin{eqnarray*}
\pi_0\map \big(\BB{\rm O}_{\nabla},\BB \ZZ/2\times \BB^2\Z/2\big)&\cong& 
\pi_0\map\big(\BB{\rm O}_{\nabla},\delta(B\ZZ/2\times B^2\Z/2)\big)
\\
&\cong & \pi_0\map\big(\Pi (\BB{\rm O}_{\nabla}),B\ZZ/2\times B^2\Z/2\big)
\\
&\cong & \pi_0\map\big(B{\rm O},B\ZZ/2\times B^2\Z/2\big)\;,
\end{eqnarray*}
and similar isomorphisms given by replacing $\BB{\rm O}_{\nabla}$ by $\BB{\rm O}$. Since the map 
$\mathcal{I}:\BB {\rm O}_{\nabla}\to \BB {\rm O}$ maps to an equivalence under $\Pi$, the claim follows. 

For the second claim, note that we have a pullback diagram of smooth stacks (without the connections)
$$
\xymatrix@C=3.8em@R=1.5em{
\BB {\rm Spin}\ar[r] \ar[d] & \ast \ar[d] 
\\
\BB {\rm O}\ar[r]^-{(w_1, \;w_2)} & \BB \ZZ/2\times \BB^2\Z/2\;.
}
$$
This follows easily from the fact that for each $n$, the double cover projection 
$p:{\rm Spin}(n)\to {\rm SO}(n)\into {\rm O}(n)$ is a smooth map and fits into the central extension of Lie 
groups $\ZZ/2\into {\rm Spin}(n)\to {\rm SO}(n)$. Moreover, since $p$ induces the isomorphism 
$\times 2: \mathfrak{so}(n)\to \mathfrak{so}(n)$ at the level of Lie algebras, one finds that 
we have a further pullback diagram 
$$
\xymatrix@C=3em@R=1.5em{
\BB {\rm Spin}_{\nabla}\ar[r] \ar[d]_-{(p\;,\;\times 2)} & \BB {\rm Spin} \ar[d]^-{p} 
\\
\BB {\rm O}_{\nabla} \ar[r] & \BB{\rm O}
}
$$
where the horizontal maps are the maps that forget the connections (i.e. they are induced by the projection 
$\Omega^1(-;\mathfrak{so}(n))\to \ast$). The result follows from the pasting law for pullbacks.
\endofproof

We now explain what orientation with respect $\widehat{\rm KO}_{\widehat{\rm tw}}$ means. 

\begin{corollary}[Differential twisted Spin structure]
A differential twisted Spin structure essentially coincides with the underlying twisted Spin structure.  
The space of $(\sigma_1,\sigma_2)$-twisted differential Spin structures is given by the homotopy fiber
$$
\xymatrix{
(w_1,w_2)\text{-}{\rm Struc}(X)_{(h_1, \;h_2)}\ar[r]\ar[d] & \ast\ar[d]^-{(\sigma_1, \;\sigma_2)}
\\
\map(X,\BB {\rm O}_{\nabla})\ar[r]^-{(w_1, \;w_2)} & \map(X,\BB \ZZ/2\times \BB^2\ZZ/2)
}
$$ 
and the bottom map $(w_1,w_2)$ is induced by the composite 
$\BB {\rm O}_{\nabla}\longrightarrow \BB{\rm O}\xrightarrow{(w_1, \;w_2)}  \BB \ZZ/2\times \BB^2\ZZ/2$ .
\end{corollary}

We now ask when such twisted structures exit. Equivalently, the question is: 
given a B-field as a cohomology class $B\in H^2(M; \Z/2)$, 
when can it satisfy anomaly cancellation? 
The following is motivated by a reinterpretation of Theorem  5.9 in \cite{KW2} (attributed to Lannes),
to  tell us when the twist is a Stiefel-Whitney class of a vector bundle.

\begin{proposition}
[Existence of twisted Spin structure and anomaly cancellation]
A necessary condition to have a twisted Spin structure and hence anomaly cancellation is 
that the B-field satisfies $\beta (B^2)=0$, where $\beta$ be the Bockstein 
$K(\Z/2, 4) \to K(\Z, 5)$, representing the cohomology operation $H^4(- ; \Z/2) \to H^5(- ; \Z)$.
 Moreover, this condition is sufficient for spacetime manifolds of 
dimension $\leq 7$. 
\end{proposition}
\theproof
We utilize the AHSS for $\op{KO}$-theory. If a degree 2 class $B\in H^2(X;\ZZ/2)$ is killed by all differentials, 
then it lifts through a map 
$$
\phi:\op{KO}^0(X)\longrightarrow H^2(X;\ZZ/2)\;.
$$ 
Since the differentials are natural with respect to pullback, standard representability arguments, along with the 
identification of the zeroth level $\op{KO}_0\simeq \op{BO}\times \ZZ$, imply that $\phi$ can be identified with the map 
\footnote{Note that the elements $H^2$ converge to elements in the filtration level $F_2\op{KO}^0(X)$. 
Such elements necessarily factor through $\op{BSO}\into \op{BO}\times \ZZ$.}
$$
\op{KO}^0(X)\cong \pi_0\map(X,\op{BO}\times \ZZ)\overset{w_2}{\longrightarrow} H^2(X;\ZZ/2)\;.
$$
Thus, both necessary and sufficient conditions for $B$ to lift through $w_2$ are provided by the differentials. 
The first nontrivial differential (at the relevant stage) in the AHSS is $d_3=\beta Sq^2$. Thus, we must have 
$$
\beta Sq^2(B)=\beta(B^2)=0\;,
$$
which gives the desired necessary condition on $B$. For manifolds of dimension $\leq 7$, all higher 
differentials vanish and hence this condition is also sufficient. 
\endofproof

\begin{remark}[Further conditions]
It is likely that the condition $\beta(B^2)=0$ is far from sufficient in general. Indeed, the higher differentials can be 
computed by considering universal examples. For $d_7$, the universal space on which $\beta Sq^2$ vanishes is 
the homotopy fiber
 $$
\xymatrix@R=1.3em{
K(\ZZ,4)\ar[r] & P\ar[r]\ar[d] & \ast\ar[d]
\\
&K(\ZZ/2,2)\ar[r]^{\beta Sq^2} &K(\ZZ,5) \;,
 }
 $$
 which is the 2-stage Postnikov section of $\op{BSO}$. If the next k-invariant in the Postnikov tower does not vanish, 
 then it can be shown that $d_7$ does not vanish. The cohomology of $P$ in degree 9 can be shown to be nonvanishing 
 (it is isomorphic to $\ZZ/2$) via the Serre spectral, and so one expects $d_7$ to be nonvanishing.
\end{remark}
%

\subsubsection*{\bf (V) Postnikov sections and type II string theory} 
 
 Twistings of \emph{complex} K-theory are classified by the set \cite{DK} 
$
 H^0(X; \Z/2) \times H^1(X; \Z/2) \times H^3(X; \Z)
$.
In \cite{DFM2}, it was shown that that these twists can be identified with the degree $-1$ cohomology of a theory $R$, 
which is represented by the Postnikov section $\op{ko} \langle 0 \cdots 4\rangle$ of $\op{ko}$, the  connective KO-theory 
spectrum. In fact, there is an exact sequence of abelian groups 
$$
\xymatrix{
0 \ar[r] & H^3(M; \Z) \ar[r] & R^{-1}(M) \ar[r]^-{(t, \; a)} & H^0(M; 
Z/2) \times H^1(M; \Z/2) \ar[r] & 0
}.
$$
The relevance of this is that the flux of the oriented superstring B-field   
is proposed to lie in $R^{-1}(M)$, which inherits both an additive and multiplicative structure from $\op{ko}$. Another novelty 
which is emphasized in \cite{DFM2} is that the twists account not just for the $B$-field, but also incorporate twists of degree 
zero and one. Twists of $\op{KO}$ theory give rise to twists of the theory $R$ and in the orientifold setting, the twists are classifies 
by the set $R^{w-1}(M)$, with $w:M\to K(\ZZ/2,1)$ classifying the orientifold double cover.  Note that analogous Postnikov sections 
have been considered beyond K-theory in \cite{ABG}. 

\medskip
Much of what is outlined in \cite{DFM1}\cite{DFM2} can, as indicated there, be extended beyond the topological 
setting and on an orientifold, the B-field can be identified as an element in the twisted differential cohomology group 
$\widehat{R}^{w-1}(X)$ for some suitable differential refinement of the cohomology theory $R$. Indeed, consider the
 differential refinement of the Postnikov section $\op{ko}\langle 0,\hdots,4\rangle$ given by the pullback 
$$
\xymatrix{
\widehat{\op{ko}}\langle 0,\hdots,4\rangle\ar[r]\ar[d] & \Omega^*(-;\tau_{\leq 4}\pi_*(\op{ko}))\ar[d]
\\
\op{ko}\langle 0,\hdots,4\rangle\ar[r] & H(\tau_{\leq 4}\pi_*(\op{ko}))
}
$$ 
where $\tau_{\leq 4}\pi_*(\op{ko})$ is the truncation of the graded ring $\pi_*(\op{ko})$ retaining only 
degrees $\leq 4$ (i.e., these are just the stable homotopy groups of $\op{ko}\langle 0,\hdots 4\rangle$). This sheaf 
of spectra defines a differential cohomology theory and we denote its value on a smooth manifold $M$ by $\widehat{R}(M)$. 
By shifting the spectrum $\op{ko}\langle 0,\hdots,4\rangle$, we perform the refinements in all degrees, 
and this gives a natural differential refinement $\widehat{R}^{-1}(X)$.  Just as $\widehat{\op{KO}}$, 
this theory can be twisted by classes in $H^1(M;\Z/2)$. Employing the explicit constructions 
of twisted differential (Deligne) cohomology in \cite{GS4}\cite{GS6}, we have the following.
 
\begin{proposition}[Refinement of the Postnikov section / spectrum of twists]
Let $X$ be an orientifold with classifying map $w:X\to \BB\ZZ/2$. Let $\mathcal{L}\to X$ be a real 
line bundle, equipped with flat connection whose corresponding monodromy representation agrees with the 
representation corresponding to $w$. We have a natural bijection of sets
$$
\widehat{R}^{w-1}(X)\cong H^0(X;\ZZ/2)\times H^1(X;\ZZ/2)\times \widehat{H}^3(X;\nabla)\:,
$$
where $\widehat{H}^3(X;\nabla)$ is degree 3 differential cohomology, twisted by the triple 
$(w,\nabla,\mathcal{L})$. 
\end{proposition}
\theproof
We apply the twisted differential AHSS. Because of the low degrees, the spectral sequence has only one 
possible nontrivial differential given by 
$$
d_3:\Omega^3(X;\mathcal{L})\longrightarrow H^3(X;\mathcal{L}/\mathcal{Z})\;.
$$
From the identification in Proposition \ref{prop-geomdiff}, this differential vanishes if and only if 
$\frac{1}{2}[\omega_4]\equiv 0 \mod \mathcal{Z}$. Therefore, we find that we have a short exact sequence
\(\label{sesrtwth}
\xymatrix{
 H^0(X;\ZZ/2)\times H^1(X;\ZZ/2)\times H^3(X;\mathcal{L}/\mathcal{Z})\ar[r] & \widehat{R}^{w-1}(X)
\ar[r] & {\rm ker}(d_3)\subset\Omega_{\rm fl}^3(X;\mathcal{L})
}.
\)
Now from \cite{GS4}, we have a short exact sequence of groups 
$$
H^3(X;\mathcal{L}/\mathcal{Z})\longrightarrow  
\widehat{H}^3(X;\nabla)\longrightarrow \mathcal{K}\subset \Omega^3_{\rm fl}(X;\mathcal{L})\;,
$$
and $\times 2:\ker(d_3)\to \mathcal{K}$ defines a bijection. Since every short exact sequence of sets splits, 
comparison with the sequence \eqref{sesrtwth} proves the claim.
\endofproof

Thus, we see that the cohomology group $\widehat{R}^{w-1}(X)$ indeed incorporates the full differential 
data of the $B$-field and not just its underlying topological class.

 \medskip
 The above was only a sample of applications and we plan more  elsewhere. 
 Our machinery can also be used in other settings, for instance,  in the context of positive scalar 
 curvature, as indicated in the Introduction.



\end{document}